\def\version{December 22, 2003}

\documentclass[12pt]{article}
\usepackage[psamsfonts]{amsfonts}
\usepackage{amsmath,amssymb,amsthm}
\usepackage[dvips]{graphicx}
\usepackage{eepic}

\def\UseSection{
        \numberwithin{equation}{section}
    \theoremstyle{plain}
        \newtheorem{theorem}    {Theorem}[section]
        \DefineTheorems 
}

\def\DefineTheorems{
    
    \newtheorem{lemma}      [theorem] {Lemma}
    
    \newtheorem{prop}       [theorem] {Proposition}
    
    \newtheorem{cor}        [theorem] {Corollary}

    \theoremstyle{definition}
    \newtheorem{defn}       [theorem] {Definition}

    \theoremstyle{definition}

}

\newcommand{\bt}   {\begin{theorem}}
\newcommand{\et}   {\end  {theorem}}
\newcommand{\bl}   {\begin{lemma}}
\newcommand{\el}   {\end  {lemma}}
\newcommand{\bp}   {\begin{prop}}
\newcommand{\ep}   {\end  {prop}}
\newcommand{\bc}   {\begin{cor}}
\newcommand{\ec}   {\end  {cor}}
\newcommand{\bd}   {\begin{defn}}
\newcommand{\ed}   {\end  {defn}}

\newcommand{\ba}   {\begin{array}}
\newcommand{\ea}   {\end  {array}}
\newcommand{\be}   {\begin{enumerate}}
\newcommand{\ee}   {\end  {enumerate}}
\newcommand{\bi}   {\begin{itemize}}
\newcommand{\ei}   {\end  {itemize}}

\def\eq#1\en{\begin{equation}#1\end{equation}}  
\def\eqsplit#1\ensplit{
    \begin{equation}\begin{split}#1\end{split}\end{equation}
    }
\def\eqalign#1\enalign{
    \begin{align}#1\end{align}
    }
\def\eqmul#1\enmul{
    \begin{multline}#1\end{multline}
    }
\newcommand{\eqarrstar} {\begin{eqnarray*}} 
\newcommand{\enarrstar} {\end{eqnarray*}} 
\newcommand{\eqarray}   {\begin{eqnarray}} 
\newcommand{\enarray}   {\end{eqnarray}} 
\newcommand{\nnb}   {\nonumber \\} 

\newcommand{\lbeq}[1]  {\label{e:#1}}
\newcommand{\refeq}[1] {\eqref{e:#1}}    

\makeatletter
\newcommand{\labelcounter}[2]{{%
    \stepcounter{#1}
    \protected@write\@auxout{}%
    {\string\newlabel{#2}{{\csname the#1\endcsname}{\thepage}}}%
    {\ref{#2}}
    }}
\makeatother
\newcommand{\sss}   { \scriptscriptstyle } 

\newcommand{\Bbold} {{\mathbb B}}

\newcommand{\Ebold} {{\mathbb E}}

\newcommand{\Pbold} {{\mathbb P}}

\newcommand{\Vbold} {{\mathbb V}}

\newcommand{\Zbold} {{\mathbb Z}}

\newcommand{\Ccal}   {\mathcal{C}} 
 
\newcommand{\Ecal}   {\mathcal{E}}

\newcommand{\Mcal}   {\mathcal{M}}

\newcommand{\Pcal}   {\mathcal{P}}

\newcommand{\Scal}   {\mathcal{S}}

\newcommand{\spose}[1] {{\hbox to 0pt{#1\hss}} }
\newcommand{\ltapprox} {\mathrel{\spose{\lower 3pt\hbox{$\mathchar"218$}}
 \raise 2.0pt\hbox{$\mathchar"13C$}}}
\newcommand{\gtapprox} {\mathrel{\spose{\lower 3pt\hbox{$\mathchar"218$}}
 \raise 2.0pt\hbox{$\mathchar"13E$}}}

\UseSection  
\renewcommand{\to}      {\rightarrow}

\newcommand{\ben}{\begin{enumerate}}
\newcommand{\een}{\end{enumerate}}

\newcommand{\Cmax} {{\Ccal}_{\rm max}}

\newcommand{\SSS}   {\sss}
\newcommand{\sM} {{\sss R}}
\newcommand{\sR} {{\sss R}}
\newcommand{\sQ} {{\sss Q}}

\newcounter{countC}  
\newcounter{countR}  
\newcommand{\sumtwo}[2]{\sum_{ \mbox{ \scriptsize
    $\begin{array}{c}
                        {#1} \\ {#2}
                        \end{array} $ }
    }
}

\newcommand{\Z}{\Zbold}

\newcommand{\conn}{\leftrightarrow}

\newcommand{\dbc}{\Leftrightarrow}
\newcommand{\ct}[1]     { \stackrel{#1}{\conn} }

\newcommand{\smallsup}[1] {{\scriptscriptstyle{({#1}})}}

\newcommand{\gr}{\mathbb G}
\newcommand{\qn}{{\mathbb Q}_n}

\newcommand{\cn}{\Omega}

\title  {
        Asymptotic expansions in $n^{-1}$ for
        \\
        percolation critical values
        on the $n$-cube and $\Z^n$
        }

\author{
Remco van der Hofstad\thanks{Department of Mathematics and Computer Science,
Eindhoven University of Technology, P.O.\ Box  513,
5600 MB Eindhoven, The Netherlands.
{\tt rhofstad@win.tue.nl}}
\and
Gordon Slade\thanks{Department of Mathematics, University of British Columbia,
Vancouver, BC V6T 1Z2, Canada. {\tt slade@math.ubc.ca}}
}

\date\version

\begin{document}

\maketitle


\begin{abstract}
We use the lace expansion to
prove that the critical values for
nearest-neighbour bond percolation on
the $n$-cube $\{0,1\}^n$ and on $\Z^n$
have asymptotic expansions, with rational coefficients, to all
orders in powers of $n^{-1}$.
\end{abstract}


\section{Main result}
\label{sec-intro}

We consider bond percolation on $\Z^n$ with edge set consisting
of pairs $\{x,y\}$ of vertices in $\Z^n$ with $\|x-y\|_1 = 1$,
where $\|w\|_1 = \sum_{j=1}^n |w_j|$ for $w \in \Z^n$.
Bonds (edges) are independently occupied with probability $p$
and vacant with probability $1-p$.
We also consider bond percolation on the $n$-cube $\qn$, which
has vertex set $\{0,1\}^n$ and edge set consisting
of pairs $\{x,y\}$ of vertices in $\{0,1\}^n$
with $\|x-y\|_1 = 1$, where we regard
$\qn$ as an additive group with addition component-wise
modulo~2.
Again
bonds are independently occupied with probability $p$
and vacant with probability $1-p$.
We write $\gr$ in place of $\qn$ and $\Z^n$ when we wish
to refer to both models simultaneously.
We write $\cn$ for the degree of $\gr$, so that
$\cn = 2n$ for $\Z^n$ and $\cn = n$ for $\qn$.

For the case of $\Z^n$, the critical value is defined by
\eq
    p_c(\Z^n) = \inf \{p : \exists \text{  an infinite
    connected cluster of occupied bonds a.s.}\}.
\en
Given a vertex $x$ of $\gr$, let $C(x)$ denote the connected
cluster of $x$, i.e., the set of vertices $y$ such that
$y$ is connected
to $x$ by a path consisting of occupied bonds.  Let $|C(x)|$
denote the cardinality of $C(x)$,
and let $\chi(p) = \Ebold_p|C(0)|$ denote the expected cluster size of
the origin.
Results of \cite{AB87,Mens86} imply that
\eq
    p_c(\Z^n) = \sup\{p : \chi(p) <\infty\}.
\en
is an equivalent definition of the critical value.

For percolation on a finite graph $\gr$, such as $\qn$, the above
characterizations of $p_c(\gr)$ are inapplicable.
In \cite{BCHSS04a,BCHSS04b,BCHSS04c} (in particular,
see \cite{BCHSS04c}),
it was shown that there is a small positive constant
$\lambda_0$ such that
the critical value $p_c(\qn) = p_c(\qn;\lambda_0)$ for the $n$-cube
is defined implicitly by
\eq
\lbeq{pcndef}
    \chi(p_c(\qn)) = \lambda_0 2^{n/3}.
\en
Given $\lambda_0$,
\refeq{pcndef} uniquely specifies
$p_c(\qn)$,
since $\chi(p)$ is a polynomial in $p$ that increases from
$\chi(0)=1$ to $\chi(1)=2^n$.  As discussed in \cite{BCHSS04c},
$p_c(\qn;\lambda_0)$
depends only weakly on the choice of $\lambda_0$.  This point is
also reflected in the Remark below.

Our main result is the following theorem.

\begin{theorem}
\label{thm-mr}
{\rm (i)}
For $\gr = \Z^n$, there are rational numbers
$a_i(\Z^n)$ such that for all $M \geq 1$,
\eq
\lbeq{aeZn}
    p_c(\Z^n) = \sum_{i=1}^M  a_i(\Z^n) (2n)^{-i}
    + O((2n)^{-M-1} )
    \quad
    \quad
    \mbox{as $n \to \infty$}.
\en
The constant in the error term depends on $M$.
\newline
{\rm (ii)}
For $\qn$, let $M \geq 1$, fix constants $c,c'$
(independent of $n$ but possibly depending on $M$), and choose $p$ such that
$\chi(p) \in [cn^{M}, c'n^{-2M}2^{n}]$.
Then there are rational numbers
$a_i(\qn)$, independent of $p,c,c'$, such that for all $M \geq 1$,
\eq
\lbeq{aeqn}
    p = \sum_{i=1}^M a_i(\qn) n^{-i} + O(n^{-M-1})
    \quad
    \quad
    \mbox{as $n \to \infty$}.
\en
The constant in the error term depends on $M,c,c'$, but does not
depend otherwise on $p$.
\end{theorem}

\noindent {\bf Remark.}
Theorem~\ref{thm-mr}(ii) states
that, for $\qn$, any $p$ for which
$\chi(p) \in [cn^M, c'n^{-2M}2^{n}]$ will have the
\emph{same} expansion
up to an error $O(n^{-M-1})$, with the error independent of $p$.
Note that the expansion \refeq{aeqn}
is valid
simultaneously to all orders $M$ if we choose
$p$ to lie eventually in all intervals
$[cn^{M}, c'n^{-2M}2^{n}]$.
Thus, \refeq{aeqn} holds
simultaneously for all $M$ for any $p$ for which $\chi(p)$ is
in the interval
$[f_n,f_n^{-1}2^{n}]$, where $f_n$ is a sequence that grows
faster than any power and more slowly than $e^{\alpha n}$ for
all $\alpha >0$ (e.g., $f_n = 2^{\sqrt{n}}$).
In particular, \refeq{aeqn} holds for $p=p_c(\qn;\lambda_0)$ for
\emph{any} constant $\lambda_0>0$, with the $a_i(\qn)$
independent of $\lambda_0$.

\medskip
The formulas \refeq{aeZn}--\refeq{aeqn} for $M=3$
were obtained in \cite{HS04a}, with $a_1(\qn)=a_1(\Z^n)=1$,
$a_2(\qn)=a_2(\Z^n)=1$, $a_3(\qn)=a_3(\Z^n)=\frac{7}{2}$,
and it was predicted (but not proved)
that $a_4(\qn)$ and $a_4(\Z^n)$ are different.
Equation~\refeq{aeZn} with $M=3$
was proved previously
in \cite{HS95,HS93up}, where it was
shown that $p_c(\Z^n) =  (2n)^{-1} + (2n)^{-2}
+ \frac{7}{2}(2n)^{-3} +O((2n)^{-4})$.

The expansion
\eq
\lbeq{pcGR}
    p_c(\Z^n)
    =  \frac{1}{2n} + \frac{1}{(2n)^2} + \frac{7}{2(2n)^3}
    +\frac{16}{(2n)^4} + \frac{103}{(2n)^5} + \cdots
\en
was reported in \cite{GR78}, but with no rigorous bound on the remainder.
The convergence of $p_c(\Z^n)$ to the leading term $(2n)^{-1}$
has been studied in \cite{ABS02,BK94,Gord91,HS90a,Kest90}, with various
error estimates.

For $\qn$, the leading term for the critical
value was identified as $n^{-1}$
in \cite{AKS82}, and this was refined in \cite{BKL92}.
It was subsequently proved (see \cite[(1.10)]{BCHSS04c})  that
\eq
\lbeq{pcqn1o}
    1-\lambda_0^{-1} 2^{-n/3}
    \leq
    np_c(\qn)
    \leq 1 +O(n^{-1}).
\en
This gives \refeq{aeqn} for $M=1$
for $p=p_c(\qn)$.

In \cite{BCHSS04a,BCHSS04b,BCHSS04c}, it was conjectured that
the phase transition for $\qn$ takes place within a scaling window
of width $2^{-n/3}$ centred at $p_c(\qn)$, and
it was proved that the scaling window
has width at most $e^{-cn^{1/3}}$.
Except in the unlikely circumstance
that the expansion \refeq{aeqn}
eventually terminates as a polynomial in $n^{-1}$,
any truncation of the expansion gives a result that lies outside this
exponentially small
scaling window for large enough $n$.  The non-perturbative definition
\refeq{pcndef} of $p_c(\qn)$ therefore tracks the scaling window
more accurately than any polynomial in $n^{-1}$ can ever do.

On the other hand, the asymptotic expansion gives information
on the phase transition on $\qn$, at every polynomial scale.
To explain this, we first recall some results from \cite{BCHSS04c}.
Let $p=p_c(\qn) + \epsilon n^{-1}$, where $\epsilon$ may depend on $n$,
and
let $\Cmax$ denote a cluster of maximal size.
By \cite[Theorem~1.1]{BCHSS04c}, if $\epsilon <0$ and $\lim_{n \to \infty}
|\epsilon| 2^{n/3} = \infty$, then
    \eqalign
    \lbeq{c1}
        \chi(p) & = \frac {1} {|\epsilon|} [1+ o(1)],
    \enalign
    \eqalign
    \lbeq{c1a}
       \frac{c_1}{\epsilon^2}
       \leq
       |\Cmax|  \leq \frac{2(\log 2)n}{\epsilon^2}  [1+o(1)]
    \quad\text{ a.a.s.},
    \enalign
where $c_1$ is a universal constant,
and where we say that a sequence $E_n$ of events holds a.a.s.\
(asymptotically almost surely) if
$\lim_{n \to \infty}\Pbold(E_n)=1$.
By \cite[Theorems~1.4, 1.5]{BCHSS04c}, there are universal constants
$c_2,c_3,c_4$ such that if $\epsilon >0$ and $\lim_{n \to \infty}
\epsilon 2^{n/3} = \infty$, then
    \eq
    \lbeq{c2}
    c_2 \epsilon^2 2^n \leq \chi(p)
    \leq c_3 \epsilon^2 2^n,
    \en
    \eq
     \lbeq{c3}
    |\Cmax|  \geq c_4 \epsilon 2^n
    \quad\text{ a.a.s.}
    \en

For $M \geq 1$, let
\eq
    p_c^\smallsup{M}(\qn) = \sum_{i=1}^M a_i(\qn) n^{-i}
\en
denote an approximate critical value.  Let $p=p_c^\smallsup{M}(\qn) +
\delta n^{-M}$, so that by \refeq{aeqn} $p=p_c(\qn) + \epsilon n^{-1}$
with $\epsilon = \delta n^{1-M} [ 1+ O(\delta^{-1} n^{-1})]$.
It follows from \refeq{c1}--\refeq{c1a} that for
fixed $\delta <0$, as $n \to \infty$, we have
    \eqalign
    \lbeq{chiasy2}
        \chi(p) & = |\delta|^{-1}n^{M-1} [1+ o(1)],
    \enalign
    \eqalign
       \frac 12 c_1 \delta^{-2} n^{2M-2}
       \leq
       |\Cmax|  \leq 2(\log 2) \delta^{-2} n^{2M-1} [1+o(1)]
    \quad\text{ a.a.s.}
    \enalign
In addition, it follows from \refeq{c2}--\refeq{c3} that for
fixed $\delta >0$,
    \eqalign
    \lbeq{chiasysup2}
     \frac 12 c_2 \delta^2 n^{2(1-M)} 2^n \leq \chi(p)
    & \leq 2  c_3 \delta^2 n^{2(1-M)} 2^n,
    \enalign
    \eqalign
     \lbeq{LCbdsup2}
    |\Cmax| & \geq  \frac 12 c_4 \delta n^{1-M} 2^n
    \quad\text{ a.a.s.}
    \enalign
This shows that
the largest cluster jumps from poly-logarithmic in the volume
for $\delta <0$, to proportional to the volume
up to poly-logarithmic factors for $\delta >0$.
Thus there is a phase transition on scale $n^{-M}$ at $p_c^\smallsup{M}(\qn)$,
for each $M \geq 1$.  A similar transition takes place if we instead
consider $p=p_c(\qn)+ \delta n^{-M}$,
where the critical point $p_c(\qn)$ is of course
{\em independent} of $M$.  Although $p_c(\qn)$
locates the phase transition
simultaneously for perturbations on scale $n^{-M}$ for all $M$,
there is nevertheless a satisfying concreteness to the polynomial
$p_c^\smallsup{M}(\qn)$ in $n^{-1}$,
compared to the implicitly defined $p_c(\qn)$ of \refeq{pcndef}.
Moreover, $p_c^\smallsup{M}(\qn)$ is independent of the parameter $\lambda_0$.

In the next proposition, we conclude
from \refeq{chiasy2} and \refeq{chiasysup2}
that Theorem~\ref{thm-mr}(ii) follows if
we can prove
\refeq{aeqn}  for any one
particular choice of $p$, such that
$\chi(p)$ lies in the interval $[cn^{M}, c'n^{-2M}2^{n}]$.
We will make a convenient choice in \refeq{pbardef} below.

\begin{prop}
\label{prop-p} Let $\gr = \qn$.
If $p$ obeys \refeq{aeqn} for some
fixed sequence $p$ (depending on $n$) such that $\chi(p)\in
[cn^{M}, c'n^{-2M}2^{n}]$, then \refeq{aeqn} is valid for every
such $p$.
\end{prop}

\proof
Define $p_1$ and $p_2$ by $\chi(p_1)=cn^M$ and $\chi(p_2)=c'n^{-2M}2^n$.
Then $p_1 \leq p_2$.  We will prove that
$p_1 \geq p_c(\qn) - O(n^{-M-1})$ and $p_2 \leq p_c(\qn) + O(n^{-M-1})$,
where the constants in the error bounds depend only on $c,c'$.
This suffices, since it implies that if $\chi(p),\chi(p') \in
[cn^{M}, c'n^{-2M}2^{n}]$, then $|p'-p| = O(n^{-M-1})$.

To prove the bound on $p_1$, we
let $\epsilon_1 = 2c^{-1}n^{-M}$.  By \refeq{c1},
for $n$ sufficiently large,
\eq
    \chi( p_1) = \frac {2}{\epsilon_1} \geq \chi(p_c(\qn)- \epsilon_1 n^{-1}),
\en
so by the monotonicity of $\chi$,
\eq
\lbeq{p1bd}
    p_1 \geq p_c(\qn)- \frac{\epsilon_1}{n} = p_c(\qn) - \frac{2}{cn^{M+1}}.
\en
For the bound on $p_2$,
we let $\epsilon_2 = \sqrt{c'/c_2}n^{-M}$.  By \refeq{c2},
\eq
    \chi( p_2) =  c_2 \epsilon_2^2 2^n
    \leq \chi(p_c(\qn)+\epsilon_2 n^{-1}),
\en
so by the monotonicity of $\chi$,
\eq
\lbeq{p2bd}
    p_2 \leq p_c(\qn)+ \frac{\epsilon_2}{n} = p_c(\qn) + \sqrt{\frac{c'}{c_2}}
    \frac{1}{n^{M+1}}.
\en
\qed

\medskip
We expect that the full asymptotic expansions
$\sum_{n=1}^\infty a_i(\gr)\cn^{-i}$ for $\qn$ and $\Z^n$ do \emph{not}
converge, although we have no serious evidence for
this conjecture. On a very naive level, the coefficients
$1,1,\frac{7}{2}, 16, 103$ stated in \refeq{pcGR} are going in the
uncomfortable direction.
Also, there is an example where divergence has been proven.
In \cite{FS90}, Fisher and Singh  review
expansions in the inverse
dimension for the critical temperature of spin systems, and
in particular for the $N$-vector model
\cite{GF74}.  It has not been proven that such expansions
are asymptotic for general $N$. However, for the spherical model \cite{BK52},
which corresponds to the limit $N \to \infty$
of the $N$-vector model, the critical temperature
$T_c(\Z^n)$ is exactly equal to
\eq
\lbeq{I10def}
        T_c(\Z^n) = \left[ \int_{[-\pi,\pi]^n}\frac{1}
        {2(n-  \sum_{j=1}^n \cos k_j)}
    \frac{dk_1}{2\pi}\cdots \frac{dk_n}{2\pi} \right]^{-1} .
\en
The integral \refeq{I10def}
has an asymptotic expansion to all orders in powers of $1/n$,
but this expansion
is divergent \cite{GF74}.   There is no reason to expect
that such divergence is limited to the spherical model.

As a side remark, we note that it is pointed out in \cite{GF74}
that there are sign changes in the expansion of $T_c(\Z^n)$
at orders $(2n)^{-12}$ and $(2n)^{-20}$.
From \refeq{pcGR}, one might guess that
$a_i(\Z^n) \geq 0$ for all $i$,
but the spherical model shows that negative coefficients
can occur relatively late in the series.  Indeed, although we do not
determine their overall signs except for $i =1,2,3$, in our proof there
are contributions to the $a_i(\gr)$ of both signs.

For the connective constant for
self-avoiding walks on $\Z^n$, existence of an asymptotic
expansion to all orders in $(2n)^{-1}$ was proved in \cite{HS95},
but the corresponding result for percolation was not obtained.
Our method is based on the lace expansion and follows the same
general approach as that used for the connective constant in
\cite{HS95}, but the details here are significantly different and
substantially more difficult.

\section{Application of the lace expansion}
\label{sec-le}

For $\qn$ or $\Z^n$ with $n$ large, the lace expansion \cite{HS90a}
gives rise to an identity
        \eq
        \lbeq{tauk'}
        \chi(p) = \frac{1+\hat \Pi_p}
        {1-\cn p[1+\hat \Pi_p]}
        ,
        \en
valid for $p \leq p_c(\gr)$.
The function $\hat\Pi_p$ is
finite for this range of $p$.  Although we do not display
the dependence explicitly in the notation, $\hat\Pi_p$ does
depend on the graph $\qn$ or $\Z^n$.
For a derivation of the lace expansion, see,
e.g., \cite[Section~3]{BCHSS04b}.
It follows from \refeq{tauk'} that
\eq
\lbeq{pcchi}
    \cn p = \frac{1}{1+\hat \Pi_{p}}
    - \chi(p)^{-1} .
\en

For $\Z^n$, \refeq{tauk'} follows from
results in \cite[Section~4.3.2]{HS90a}.
(Note the notational difference that
in \cite{HS90a} what we are calling here
$\hat{\Pi}_p$ is equal to $\sum_{n=0}^\infty (-1)^n \hat{g}_n(0)$.)
Since $\chi(p_c(\Z^n))=\infty$, \refeq{pcchi} gives
\eq
\lbeq{pcdefZd}
    2np_c(\Z^n) =  \frac{1}{1+\hat \Pi_{p_c(\Z^n)}}.
\en
Bounds of \cite{HS90a} imply that $|\hat \Pi_{p_c(\Z^n)}|=O(n^{-1})$.
This gives
$p_c(\Z^n) = (2n)^{-1} + O(n^{-2})$, which is \refeq{aeZn} for $M=1$.

For $\qn$, the identity
\refeq{tauk'} is established in \cite[(6.1)]{BCHSS04b}.
We fix a sequence $f_n$ as in the remark below Theorem~\ref{thm-mr}.
That is, we require that $\lim_{n\to \infty}f_nn^{-M} =\infty$
for every positive integer $M$ and that $\lim_{n\to\infty}f_ne^{-\alpha n}
=0$ for every $\alpha >0$.
In view of Proposition~\ref{prop-p}, it suffices to prove
\refeq{aeqn} for $p=\bar p$, where $\bar p$ is defined by
\eq
\lbeq{pbardef}
    \chi(\bar p)=f_n.
\en
For $p = \bar p$,
\refeq{pcchi} gives
    \eq
    \lbeq{pcform}
    n \bar p
    = \frac{1}{1+\hat \Pi_{\bar p}}
    +O(f_n^{-1}) .
    \en
The second term on the right hand side is an error term for
\refeq{aeqn}, and can be neglected in the proof of
Theorem~\ref{thm-mr}(ii).
It follows from results in \cite{BCHSS04b} that
$|\hat\Pi_{\bar p}| \leq O(n^{-1})$.
In more detail, it follows from \cite[Proposition~5.2]{BCHSS04b}
that $|\hat{\Pi}_p| \leq
\mbox{const}(\lambda^3 \vee \beta)$, where
$\lambda = \chi(p) 2^{-n/3} \leq f_n 2^{-n/3}$ for
$p \leq \bar p_c(\qn)$, and $\beta$ is proportional to $n^{-1}$
by \cite[Proposition~2.1]{BCHSS04b}.
With \refeq{pcform}, this implies that $\bar p = n^{-1} +O(n^{-2})$,
which is \refeq{aeqn} for $M=1$.

Henceforth, we will write
\eq
    \bar p_c = \bar p_c(\gr)=
    \begin{cases}
    \bar p & (\gr = \qn)
    \\
    p_c(\Z^n) & (\gr = \Z^n).
    \end{cases}
\en
The identities \refeq{pcdefZd} and \refeq{pcform}  give recursive
equations for $\bar p_c$, in which an input for $\bar p_c$ on the right hand
side gives rise to an improved value of $\bar p_c$ on the left hand side.
To prove Theorem~\ref{thm-mr} using this recursion, we will
apply the following proposition.

\begin{prop}
\label{prop-Piae}
Fix $M \geq 1$.
For $\gr = \qn$ and $\gr = \Z^n$,
there are rational numbers $\alpha_{j,i,M}=\alpha_{j,i,M}(\gr)$ and a positive
integer $L_M$ such that
for $p \leq \bar p_c(\gr)$,
\eq
\lbeq{Piasy}
     \hat \Pi_{p}
    =
    \sum_{i=1}^{L_M} \sum_{j=0}^{i-1} \alpha_{j,i,M} \cn^j
    p^i
    +
    O(\cn^{-M-1}).
\en
The constant in the error term depends on $M$.
\end{prop}

\noindent {\em Proof of Theorem~\ref{thm-mr} assuming Proposition~\ref{prop-Piae}.}
We restrict attention to $p = \bar p_c$, which is sufficient by
Proposition~\ref{prop-p}.
The proof is by induction on $M$.  As discussed above,
we know that \refeq{aeZn} and \refeq{aeqn}
hold for $M=1$.  We assume that \refeq{aeZn}--\refeq{aeqn} hold
for $1,\ldots, M$, and prove
the corresponding result for $M+1$.

By \refeq{pcdefZd} and \refeq{pcform},
\eq
\lbeq{pcM}
    \bar p_c =
    \frac{1}{\cn}
    \frac{1}{1+\hat\Pi_{\bar p_c}}
    + O(\cn^{-M-2})
\en
(in fact the error term is zero for $\Z^n$ and is
$O(n^{-1}f_n^{-1})$ for $\qn$).
By the induction hypothesis,
there are rational $\beta_{k,i}=\beta_{k,i}(\gr)$ such that
for $i = 1,\ldots , L_M$,
\eqalign
    \bar p_c^i
    & = \left[ \sum_{j=1}^M a_j \cn^{-j} + O(\cn^{-M-1})\right]^i
    = \cn^{-i}
    \left[ \sum_{k=0}^{M-1} \beta_{k,i} \cn^{-k} + O(\cn^{-M})\right]
    .
\enalign
By Proposition~\ref{prop-Piae}, this implies that
\eqalign
    \hat \Pi_{\bar p_c}
    & =
     \sum_{i=1}^{L_M} \sum_{j=0}^{i-1} \alpha_{j,i,M} \cn^j
    \cn^{-i}
    \left[ \sum_{k=0}^{M-1} \beta_{k,i} \cn^{-k} + O(\cn^{-M})\right]
    +
    O(\cn^{-M-1})
    \nnb
    \lbeq{PiM1}
    & =
    \sum_{l=1}^{M}  \gamma_{l,M} \cn^{-l}
    +
    O(\cn^{-M-1}),
\enalign
for some rational coefficients $\gamma_{l,M}=\gamma_{l,M}(\gr)$.
Substitution of \refeq{PiM1} into \refeq{pcM} shows that there
are rational coefficients $a_{j,M}$ such that
\eq
\lbeq{pcasypf}
    \bar p_c = \sum_{j=1}^{M+1} a_{j,M} \cn^{-j} + O(\cn^{-M-2}).
\en
The consistency of \refeq{pcasypf} with the induction hypothesis
implies that for $j=1,\ldots, M$ the coefficients $a_{j,M}$ in fact
depend only on $j$.  This completes the proof.
\qed

Note that the precise value of $L_M$ is not needed in the above proof.
It remains to prove Proposition~\ref{prop-Piae}.
For this, we will use the description of $\hat\Pi_p$ given
in the next section.

\section{The function $\hat\Pi_p$}
\label{sec-Pi}

The function $\hat\Pi_p$ has the form
    \eq
    \lbeq{2pt.37}
    \hat\Pi_p = \sum_{N=0}^{\infty} (-1)^N
    \hat\Pi_p^{\smallsup{N}}.
    \en
To define $\hat\Pi_p^{\smallsup{N}}$
(see \cite{BCHSS04b} or \cite{HS90a}), we need the following definition.

   \begin{defn}
    \label{def-inon}
    (i)
Given a bond configuration, we say that $x$ is {\em doubly connected
    to}\/ $y$, and write $x \dbc y$, if $x=y$ or
    if there are at least two bond-disjoint
    paths from $x$ to $y$
    consisting of occupied bonds.
    \newline
        (ii)
        Given a bond configuration, vertices $x,y$,
        and a set $A$ of vertices of $\gr$, we
        say $x$ and $y$ are \emph{connected through $A$}, and write
        $x \ct{A} y$, if every
        occupied path connecting $x$ to $y$ has at least one bond
        with an endpoint in $A$.
    \newline
        (iii)
        Given a bond configuration, and a bond $b$, we define
        $\tilde{C}^{b}(x)$ to be the set of vertices connected to $x$
        in the new configuration obtained by setting $b$ to be vacant.
        \newline
    (iv) Given a bond configuration and vertices $x,y$,
    we say that the directed bond $(u,v)$ is {\em pivotal} for $x \conn y$
    if (a) $x \conn y$ occurs when the bond $\{u,v\}$ is set occupied, and
    (b) $x \conn y$ does not occur when $\{u,v\}$ is set vacant,
    but $x\conn u$ and $v \conn y$ do occur. (Note that there is a distinction
    between the events $\{(u,v)$ is pivotal for $x \conn y\}$
    and $\{(u,v)$ is pivotal for $y \conn x\}
    = \{(v,u)$ is pivotal for $x \conn y\}$.)
\newline (v)
Given vertices $v,x$ and a set of vertices $A$, let
    \eqalign
    \lbeq{317}
    E'(v, x; A) & = \{ v \ct{A} x\} \cap
    \{\not\exists \text{ occupied pivotal $(u', v')$ for $v \conn x$ s.t.
    $v \ct{A} u'$} \}.
    \enalign
We refer to the second event on the right hand side of \refeq{317} as
the ``NP'' (no pivotal) condition.
    \end{defn}

By definition,
\eq
\lbeq{pi0defz}
    \hat\Pi^\smallsup{0}_p = \sum_{x \neq 0} \Pbold_p (0 \dbc x)
\en
and
    \eq
    \lbeq{Pi1defa}
    \hat\Pi_p^\smallsup{1}=\sum_{x}\sum_{(u_0, v_0)}
    p\Ebold_{\sss 0} \left[ I[0 \dbc u_0]
    {\mathbb E}_{\sss 1} I[E'(v_0, x;\tilde{C}^{(u_0,v_0)}_{\SSS 0}(0))] \right].
    \en
On the right hand side
of \refeq{Pi1defa}, the cluster $\tilde{C}^{(u,
v)}_{\SSS 0}(0)$ is \emph{random} with respect to the expectation
$\Ebold_{\SSS 0}$, but
$\tilde{C}^{(u,v)}_{\SSS 0}(0)$ should be regarded as a \emph{fixed} set
inside the probability $\Pbold_{\SSS 1}$.
The latter introduces a second percolation model which
is dependent on the original percolation model only via the set
$\tilde{C}^{(u, v)}_{\SSS 0}(0)$.

In general, for $N \geq 1$,
    \eqalign
    \lbeq{2pt.38}
    \hat\Pi^{\smallsup{N}}_p
    & = 
    \sum_x \sum_{(u_0, v_0)}\cdots \sum_{(u_{N-1},v_{N-1})}
    p^N
    \Ebold_{\sss 0} I[0 \dbc u_{0}] \, \,
    \\ \nonumber &
    \quad \quad \times
    \Ebold_{\sss 1} I[E'(v_{0}, u_{1}; \tilde{C}_{0})]
     \cdots
     \Ebold_{\sss N-1} I[E'(v_{N-2}, u_{N-1}; \tilde{C}_{N-2})]
    \Ebold_{\sss N} I[E'(v_{N-1}, x; \tilde{C}_{N-1})],
    \enalign
where we have used the abbreviation
$\tilde{C}_{j} = \tilde{C}_j^{(u_{j}, v_{j})}(v_{j-1})$
(with $v_{-1}=0$), and where
each sum over $(u,v)$ is a sum over all directed bonds.
The expectations in \refeq{2pt.38} are mutually dependent
through the $\tilde{C}$ clusters.
We use subscripts for
$\tilde{C}$ and the expectations,
to indicate to which expectation $\tilde{C}$ belongs,
and refer to the bond configuration corresponding to expectation~$j$
as the ``level-$j$'' configuration.
We also write $F_j$ to indicate an event $F$ at level-$j$.

We will rewrite \refeq{2pt.38} as follows.
We set $u_N=x$ and $v_{-1}=0$, and write
\eq
    E_0 = \{0 \dbc u_0\}_0,
\en
\eq
    E_j = E'(v_{j-1},u_j,\tilde{C}_{j-1})_j
    \quad \quad (j = 1,\ldots, N),
\en
and
\eq
\lbeq{EN-def}
    E^\smallsup{N} = E_0 \cap E_1 \cap \cdots \cap E_N.
\en
We drop the
subscript $p$ on $\hat\Pi^\smallsup{N}$, and write
$\Ebold^\smallsup{N}$ for the joint expectation
$\Ebold_{\sss 0}\Ebold_{\sss 1}\cdots \Ebold_{\sss N}$.
We write the sum over $x$ and all $(u_j,v_j)$ ($j=0,\ldots,N-1)$ as
$\sum_{x,(u_j,v_j)}$.  With this notation, \refeq{2pt.38} takes the
more compact form
\eq
\lbeq{Picomp}
    \hat\Pi^\smallsup{N}_p = \sum_{x,(u_j,v_j)}
    p^N  \Ebold^\smallsup{N}
    \left[ I[E^\smallsup{N}] \right].
\en

It is known that for all $N \geq 0$,
    \eq
    \lbeq{Pibds}
    0 \leq \hat \Pi_p^{\smallsup{N}}
    \leq \left(\frac{C}{\cn }\right)^{N\vee 1}
    \quad
    \text{uniformly in $p \leq \bar p_c(\gr)$}.
    \en
For $\qn$, \refeq{Pibds} is established
in \cite[Lemma~5.4]{BCHSS04b}.  In more detail, \cite[Lemma~5.4]{BCHSS04b}
states that $\hat{\Pi}_p^\smallsup{N} \leq
[\mbox{const}(\lambda^3 \vee \beta)]^{N\vee 1}$, where
$\lambda = \chi(p) 2^{-n/3} \leq cf_n 2^{-n/3}$ for
$p \leq \bar p_c(\qn)$.
In addition,
it is shown in \cite[Proposition~2.1]{BCHSS04b} that
$\beta$ can be chosen proportional to $n^{-1}$.  This gives
\refeq{Pibds} for $\qn$.
For $\Z^n$, \refeq{Pibds} follows from
results in \cite[Section~4.3.2]{HS90a} (we emphasize again
that there are notational differences in \cite{HS90a}).

It follows immediately from \refeq{Pibds} that
to prove \refeq{Piasy} it is sufficient to prove
that there are rational numbers $\alpha_{j,i,M}$ such that
\eq
\lbeq{Piasybis}
    \sum_{N=0}^M (-1)^N \hat \Pi^\smallsup{N}_p
    =
    \sum_{i=1}^{L_M}\sum_{j=0}^{i-1}\alpha_{j,i,M} \cn^j
    p^i + O(\cn^{-M-1}).
\en
Thus it is sufficient to show that for each $N \leq M$ there
are rational numbers $\alpha_{j,i,M}^\smallsup{N}$ such that
\eq
\lbeq{PiNidentity}
    \hat \Pi^\smallsup{N}_p
    =
    \sum_{i=1}^{L_M}\sum_{j=0}^{i-1}\alpha_{j,i,M}^\smallsup{N} \cn^j
    p^i + O(\cn^{-M-1}).
\en
Suppose that we could show instead that
\eq
\lbeq{PiNasy}
    \hat \Pi^\smallsup{N}_p
    =
    \sum_{i=0}^{L_M}\sum_{j=0}^{L_M}
    \alpha_{j,i,M}^\smallsup{N} \cn^j
    p^i
    +
    O(\cn^{-M-1}).
\en
The sum on the right hand side can be rewritten as
\eq
    \sum_{k=-L_M}^{L_M} \cn^k \left( \sum_{i=-k \vee 0}^{L_M-k}
    \alpha_{i+k,i,M}^\smallsup{N} (\cn p)^i \right).
\en
It follows from \refeq{Pibds} that
\eq
     \sum_{i= 0}^{L_M-k}
    \alpha_{i+k,i,M}^\smallsup{N} (\cn p)^i
    =0
    \quad \quad
    (k >-(N \vee 1)),
\en
for a continuous range of $p$ values of order $\cn^{-1}$.
This then implies, in particular,
that $\alpha_{i+k,i,M}^\smallsup{N}=0$ for $k\geq 0$.
Thus, \refeq{PiNasy} implies \refeq{PiNidentity}.
The remainder of the paper is devoted to the proof of \refeq{PiNasy},
for fixed $N \leq M$.

In estimating error terms, it is convenient to work with
an upper bound for the event $E^\smallsup{N}$ defined in
\refeq{EN-def}. Let $E\circ F$ denote disjoint occurrence of the events
$E$ and $F$.
We define
    \eqalign
    F_0(0,u_0,w_0,z_1) &= \{0\conn u_0\}\circ \{0\conn w_0\}\circ \{w_0\conn u_0\}
    \circ \{w_0\conn z_1\},
    \lbeq{Fdefa}\\
    \!\!\!F'(v_{i-1},t_i,z_i,u_i, w_i, z_{i+1})
    &= \{v_{i-1}\conn t_i\}\circ \{t_i\conn z_i\}\circ \{t_i\conn w_i\}\nonumber\\
    &\qquad \circ \{z_i\conn u_i\}
    \circ \{w_i\conn u_i\} \circ \{w_i\conn z_{i+1}\},
    \lbeq{F1defa}\\
    \!\!\!F''(v_{i-1},t_i,z_i,u_i, w_i, z_{i+1})
    &= \{v_{i-1}\conn w_i\}\circ \{w_i\conn t_i\}\circ \{t_i\conn z_i\}\nonumber\\
    &\qquad\circ \{t_i\conn u_i\} \circ \{z_i\conn u_i\} \circ \{w_i\conn z_{i+1}\},
    \lbeq{F2defa}\\
    \!\!\!F(v_{i-1},t_i,z_i,u_i, w_i, z_{i+1})
    &=F'(v_{i-1},t_i,z_i,u_i, w_i, z_{i+1})
    \cup F''(v_{i-1},t_i,z_i,u_i, w_i, z_{i+1}),
    \\
     F_N(v_{N-1},t_N,z_N,x)
     &= \{v_{N-1}\conn t_N\}\circ \{t_N\conn z_N\}\circ \{t_N\conn x\}
    \circ \{z_N\conn x\}.
    \lbeq{Fdefsa}
    \enalign
We also define $F^\smallsup{0} = \{0\dbc u_0\}$,
and, for $N \geq 1$, let
\eq
\lbeq{FN-def}
    F^\smallsup{N} =
    \bigcup_{\vec{t},\vec{w}, \vec{z}}
    \Big(
    F_0(0,u_0,w_0,z_1)_0
    \cap \big( \bigcap _{i=1}^{N-1} F(v_{i-1},t_i,z_i,u_i, w_i, z_{i+1})_i
    \big) \cap
    F_N(v_{N-1},t_N,z_N,x)_N\Big),
\en
where $\vec{t}= (t_1,\ldots,t_N)$, $\vec{w}= (w_0,\ldots,w_{N-1})$
and $\vec{z}= (z_1,\ldots,z_N)$.
It is then a standard estimate (see \cite[Section~4.1]{BCHSS04b} for a discussion
using our present notation) that for $N \geq 0$,
    \eqalign
    \lbeq{Fbda}
    E^\smallsup{N}
    & \subset  F^\smallsup{N}.
    \enalign

\section{The polynomial}
\label{sec-poly}

In this section, we prove \refeq{PiNasy}
(and hence Proposition~\ref{prop-Piae})
apart from the estimation
of several error terms that are identified below.  These error terms are
bounded in Section~\ref{sec-error}, to complete the proof of
\refeq{PiNasy}.  We begin with some definitions.

\subsection{Basic concepts}

For $N \geq 0$, let $\gr^\smallsup{N}=(\gr(0),\ldots,\gr(N))$
denote $(N+1)$ copies of $\gr$. This graph contains the edge set
of $\gr$ within each component $\gr(i)$, and no edge links any vertex of
$\gr(i)$ with any vertex of $\gr(j)$ when $i \neq j$.
We consider percolation on $\gr^\smallsup{N}$, which is the
natural setting for the event $E^\smallsup{N}$ of \refeq{EN-def}.
We denote the set of vertices of $\gr^\smallsup{N}$ by
$(\Vbold(0),\ldots,\Vbold(N))$, where $\Vbold(j)$ is the vertex set
of $\gr(j)$.
In the event $E_j=E'(v_{j-1},u_j;\tilde{C}_{j-1})_j$, the random set
$\tilde{C}_{j-1} \subset \Vbold(j-1)$ is identified with the
corresponding subset of $\Vbold(j)$, in \refeq{317}.

\begin{defn}
\label{def-on}
Let $N \geq 0$.
Given $A=(A(0),\ldots,A(N))$ with each $A(j)\subset \Vbold(j)$,
and an event $E$
on $\gr^\smallsup{N}$,
the event $\{E$ {\em on} $A\}$
is the set of configurations for which $E$ occurs in the
possibly modified configuration in which each bond with one
or more endpoints not in $A$ is set vacant.
\end{defn}

We write $A^c=(A(0)^c,\ldots,A(N)^c)$, and, for $B=(B(0),\ldots,B(N))$,
we define $A \cap B = (A(0)\cap B(0),\ldots,A(N) \cap B(N))$.
It is a direct consequence of Definition~\ref{def-on} that
\eqalign
\lbeq{on/inprop1}
    \{E^c \text{ on } A\} &= \{E \text{ on } A\}^c,
    \\
\lbeq{on/inprop2}
    \{E \cap F \text{ on } A\} &= \{E \text{ on } A\}
    \cap \{F \text{ on } A\},
    \\
\lbeq{on/inprop3}
    \{\{E  \text{ on } A\} \text{ on }  B\} &= \{E \text{ on } A \cap B\},
\enalign
so the notion ``on $A$'' is well behaved with respect to the
operations of set theory.
In addition, the proof of \refeq{Fbda} also shows that
\eq
\lbeq{Fbdon}
    \{E^\smallsup{N} \text{ on } A\}
    \subset
    \{F^\smallsup{N} \text{ on } A\}
    \subset
    F^\smallsup{N}.
\en

\begin{defn}
\label{def-bb}
Fix a bond configuration on $\gr^\smallsup{N}$, fix the summation
variables $x$, $(u_j,v_j)$ ($j=0,\ldots,N-1$) of \refeq{Picomp},
and let $u_N=x$, $v_{-1}=0$.
\newline
(i)
The {\em backbone} at level-$j$ is the random set of vertices
defined by
\eq
    B_j = \left\{ y \in \Vbold(j)
    : \{v_{j-1} \conn y\}_j \circ \{ y \conn u_j\}_j \right\}
    \quad
    \quad
    (j=0,\ldots,N).
\en
The {\em extended backbone} $B_j^+$ at level-$j$
is
\eq
    B_j^+ =
    \begin{cases}
    B_j \cup \{y \in \tilde{C}_{j}:
    \{v_{j-1} \conn y\}_j \circ \{y \conn B_{j+1}\}_j\}
     &
    (j=0,\ldots,N-1)
    \\
    B_N & (j=N).
    \end{cases}
\en
(ii)
The {\em dimension}
$\vec D = (D_1,\ldots, D_n)\in \{0,1\}^n$ is defined by
setting $D_i=1$ if
there is a $j\in \{0,\ldots,N\}$ and a $y \in B_j^+$ with $y_i\neq 0$,
and otherwise $D_i=0$.
Let ${\cal P}_j$ denote the collection of paths consisting
of level-$j$ paths from $v_{j-1}$ to $u_j$, and level-$j$ paths
in $\tilde{C}_j$ from $v_{j-1}$ to $B_{j+1}$.
Then $\vec D$ indicates all directions in $\gr$ that
are explored by occupied paths in $\Pcal_j$,
for all $j$.
\end{defn}

Given $R>0$,
let $\Bbold_\sM = \{x \in \Vbold : \|x\|_\infty \leq R\}$.
Let $\Bbold_\sM^\smallsup{N}=(\Bbold_\sM(0),\ldots,\Bbold_\sM(N))$
denote $(N+1)$
copies of $\Bbold_\sM$.
Given $\vec d \in \{0,1\}^n$, let
\eqalign
    \Vbold_{\vec d} & = \{ x \in \Vbold : x_i=0 \text{ if } d_i=0\},
    \\
    \Vbold_{\vec d}^\smallsup{N}
    & = (\Vbold_{\vec d}(0),\ldots,\Vbold_{\vec d}(N))
    =\text{$(N+1)$ copies of $\Vbold_{\vec d}$},
    \\
    \Vbold_{\vec d,\sM}^\smallsup{N}
    & = \Vbold_{\vec d}^\smallsup{N} \cap \Bbold_{\sM}^\smallsup{N}.
\enalign
Taking $\vec d$ to be the random vector $\vec D$, we extend
Definition~\ref{def-on} to the case where $A$ is the
{\em random} set $\Vbold_{\vec D,\sM}^\smallsup{N}$
by the disjoint union
\eq
\lbeq{EonV}
    \{E \text{ on }\Vbold_{\vec D,\sM}^\smallsup{N}\}
    =
    \bigcup_{\vec d \in \{0,1\}^n}
    \big(\{E \text{ on }\Vbold_{\vec d,\sM}^\smallsup{N}\}
    \cap \{\vec D = \vec d\} \big).
\en

We will use the following lemma.

\begin{lemma}
\label{lem-onG}
For $N\geq 0$ and $R>0$,
\eq
\lbeq{onclaima}
    \{E^\smallsup{N} \text{ on } \Bbold_{\sM}^\smallsup{N}\}
    =
    \{E^\smallsup{N} \text{ on } \Vbold_{\vec D,\sM}^\smallsup{N}\}
    .
\en
\end{lemma}

\proof
The basic step in the proof is to observe
that the event $E^\smallsup{N}$
is determined by occupied paths in the extended backbones $B_j^+$
$(j=0,1,\ldots,N)$, and by the definition of $\vec D$, these paths lie
in $\Vbold_{\vec D}^\smallsup{N}$.  From this, we see that
\eq
\lbeq{Ed}
    E^\smallsup{N} \cap \{\vec D = \vec d\}
    =
    \{E^\smallsup{N} \text{ on } \Vbold_{\vec d}^\smallsup{N}\}
    \cap \{\vec D = \vec d\}.
\en
Therefore,
\eqalign
\lbeq{onclaimz}
    \{E^\smallsup{N} \text{ on } \Bbold_{\sM}^\smallsup{N}\}
    & =
   \bigcup_{\vec d \in \{0,1\}^n}
   \{E^\smallsup{N} \text{ on } \Bbold_{\sM}^\smallsup{N}\}
   \cap \{\vec D = \vec d\}
    \nnb
    & =
    \bigcup_{\vec d \in \{0,1\}^n}
    \{\{E^\smallsup{N} \text{ on } \Vbold_{\vec d}^\smallsup{N}\}  \text{ on }
    \Bbold_{\sM}^\smallsup{N}\}
    \cap \{\vec D = \vec d\}
    \nnb
    & =
    \bigcup_{\vec d \in \{0,1\}^n}
    \{E^\smallsup{N} \text{ on } \Vbold_{\vec d,\sM}^\smallsup{N}\}
    \cap \{\vec D = \vec d\}
    \nnb
    & =
    \{E^\smallsup{N} \text{ on } \Vbold_{\vec D,\sM}^\smallsup{N}\},
\enalign
where we used \refeq{Ed} for the second equality, \refeq{on/inprop3}
for the third, and \refeq{EonV} for the fourth.
\qed

\subsection{The iteration}

Henceforth, we write $\hat\Pi_p^\smallsup{N}$ as
$\hat\Pi^\smallsup{N}$, and we will introduce new subscripts with
different meaning than $p$.  Recall the definition of $\hat\Pi^\smallsup{N}$
in \refeq{Picomp}.
For $R_0>0$, we write
\eqalign
    \hat\Pi^\smallsup{N} &= \hat\Pi^\smallsup{N}_{\sM_0}
    +
    \Ecal_{1,\sM_0}^\smallsup{N},
\enalign
where
\eqalign
\lbeq{PiM}
    \hat\Pi^\smallsup{N}_{\sM_0} &= \sum_{x,(u_j,v_j)}
    p^N  \Ebold^\smallsup{N}
    \left[ I[E^\smallsup{N} \text{ on } \Bbold_{\sM_0}^\smallsup{N}]
    \right],
\\
\lbeq{E1adef}
    \Ecal_{1,\sM_0}^\smallsup{N} &= \sum_{x,(u_j,v_j)}
    p^N  \Ebold^\smallsup{N}
    \left[ I[E^\smallsup{N}]
    - I[E^\smallsup{N} \text{ on } \Bbold_{\sM_0}^\smallsup{N}] \right].
\enalign
We will choose $R_0=R_0(M)$ in Lemma~\ref{lem-E1} below,
in such a manner that
\eq
\lbeq{E1a}
    \Ecal_{1,\sM_0}^\smallsup{N} = O(\cn^{-M-1}),
\en
so this term is an error term.  The bound \refeq{E1a} is a reflection
of the fact that configurations that leave $\Bbold_{\sM_0}^\smallsup{N}$,
with $R_0$ large depending on $M$,
must contain a long extended backbone path,
and this gives rise to an error term.  This notion will be formalized
in Proposition~\ref{prop-PibdsM} below.

By Lemma~\ref{lem-onG},
\eq
\lbeq{PiMc}
    \hat\Pi^\smallsup{N}_{\sM_0} =
    \sum_{x,(u_j,v_j)}
    p^N  \Ebold^\smallsup{N}
    \left[ I[E^\smallsup{N} \text{ on } \Vbold_{\vec D,\sM_0}^\smallsup{N}]
    \right] .
\en
Let $\|\vec d\| = \sum_{i=1}^n|d_i|$ denote the $\ell^1$ norm.
Given $R_0,R_1>0$, we write
\eq
    \hat\Pi^\smallsup{N}_{\sM_0} =
    \hat\Pi^\smallsup{N}_{\sM_0, \sM_1} +\Ecal_{2,\sM_0, \sM_1}^\smallsup{N},
\en
where
\eqalign
    \hat\Pi^\smallsup{N}_{\sM_0, \sM_1}
    & =
    \sum_{x,(u_j,v_j)}
    p^N  \Ebold^\smallsup{N}
    \left[ I[E^\smallsup{N} \text{ on } \Vbold_{\vec D,\sM_0}^\smallsup{N}]
    I[\|\vec D\| \leq R_1] \right],
    \\
    \lbeq{E2def}
    \Ecal_{2,\sM_0, \sM_1}^\smallsup{N}
    & =
    \sum_{x,(u_j,v_j)}
    p^N  \Ebold^\smallsup{N}
    \left[ I[E^\smallsup{N} \text{ on } \Vbold_{\vec D,\sM_0}^\smallsup{N}]
    I[\|\vec D\| > R_1] \right].
\enalign
We show below in Corollary~\ref{cor-T2}
that for $R_0=R_0(M)$ chosen as above,
and for suitably chosen $R_1=R_1(M)$,
\eq
\lbeq{E2a}
    \Ecal_{2,\sM_0, \sM_1}^\smallsup{N} = O(\cn^{-M-1}),
\en
so this term is an error term.  The bound \refeq{E2a}
will follow from the fact that large $\|\vec D\|$ implies
the existence either
of a long extended backbone path, or of many distinct extended
backbone paths.

By definition,
\eq
\lbeq{PiEE}
    \hat\Pi^\smallsup{N}
    =
    \hat\Pi^\smallsup{N}_{\sM_0, \sM_1}
    + \Ecal_{1,\sM_0}^\smallsup{N}
    + \Ecal_{2,\sM_0, \sM_1}^\smallsup{N},
\en
with
\eqalign
\lbeq{PiMd}
    \hat\Pi^\smallsup{N}_{\sM_0, \sM_1}
    & =
    \sum_{\vec d_1 :\|\vec d_1\| \leq R_1}
    \sum_{x,(u_j,v_j)}
    p^N  \Ebold^\smallsup{N}
    \left[ I[E^\smallsup{N} \text{ on } \Vbold_{\vec d_1,\sM_0}^\smallsup{N}]
    I[\vec D = \vec d_1] \right].
\enalign
Our goal is to rewrite the expectation in
\refeq{PiMd}, up to error terms, as an expectation of an event
that occurs on some $\Vbold_{\vec d, \sM}$ with $\|\vec d\|$ and $R$
bounded depending only on $M$.  The factor $I[\vec D = \vec d_1]$ is
not yet of this form, and needs to be rewritten.

We write $\vec d > \vec d'$ if $d_i > d_i'$ for all $i =1,\ldots, n$.
Given $\bar Q$ (large), $\vec R = (R_0,\ldots,R_{\bar Q +1})$ with
$R_0 < \cdots <R_{\bar Q+1} $,
$q \in \{1,\ldots, \bar Q\}$,
and $\|\vec d_q\| \leq R_q$, we proceed as
follows.
First, we make the decomposition
\eqalign
    I[\vec D =\vec d_q]
    &  =
    I[\{\vec D = \vec d_q \}
    \text{ on } \Vbold_{\vec d_q,\sM_{q-1}}^\smallsup{N}]
    I[\vec D =\vec d_q]
\lbeq{ondecomp}
    +
    I[ \{\vec D < \vec d_q \}
    \text{ on } \Vbold_{\vec d_q,\sM_{q-1}}^\smallsup{N}]
    I[\vec D =\vec d_q],
\enalign
In the first term on the
right hand side of \refeq{ondecomp}, we make the replacement
\eq
\lbeq{Dgd}
    I[\vec D =\vec d_q] = I[\vec D \geq \vec d_q] - I[\vec D > \vec d_q].
\en
Since $\{\{\vec D = \vec d_q \}$ on
$\Vbold_{\vec d_q,\sM_{q-1}}^\smallsup{N} \} \subset
\{\vec D  \geq \vec d_q\}$, we obtain
\eqalign
    I[\vec D =\vec d_q]
    &  =
    I[\{\vec D = \vec d_q \}
    \text{ on } \Vbold_{\vec d_q,\sM_{q-1}}^\smallsup{N}]
    \lbeq{ondecomp2}
    \nnb
    & \quad +
    I[ \{\vec D < \vec d_q \}
    \text{ on } \Vbold_{\vec d_q,\sM_{q-1}}^\smallsup{N}]
    I[\vec D =\vec d_q]
    \nnb
    & \quad -
    I[\{\vec D = \vec d_q \}
    \text{ on } \Vbold_{\vec d_q,\sM_{q-1}}^\smallsup{N}]I[\vec D >\vec d_q],
\enalign
In the last term on the right hand side,
we insert the factor
\eq
    1 = I[\|\vec D\| > R_{q+1}] +  I[\|\vec D\| \leq R_{q+1}].
\en
This gives
\eqalign
    I[\vec D = \vec d_q] & =
    I[\{\vec D = \vec d_q\} \text{ on } \Vbold_{\vec d_q, R_{q-1}}^\smallsup{N}]
    \nnb & \quad
    +
    I[\{\vec D < \vec d_q\} \text{ on } \Vbold_{\vec d_q, R_{q-1}}^\smallsup{N}]
    I[\vec D = \vec d_q]
    \nnb
    &  \quad -
    I[\{\vec D = \vec d_q\} \text{ on } \Vbold_{\vec d_q, R_{q-1}}^\smallsup{N}]
    I[\|\vec D\| > R_{q+1}]
    \nnb & \quad
\lbeq{Iit}
    -
    \sumtwo{\vec d_{q+1} : \vec d_{q+1}>\vec d_q}{ \|\vec d_{q+1}\| \leq R_{q+1}}
    I[\{\vec D = \vec d_q\} \text{ on } \Vbold_{\vec d_q, R_{q-1}}^\smallsup{N}]
    I[\vec D = \vec d_{q+1}],
\enalign
since in the third term
$\vec D > \vec d_q$ follows from the facts that
$\|\vec D\| > R_{q+1}>R_q\geq \|\vec d_q\|$ and $\vec D \geq \vec d_q$.

We define $\Scal_{ \sQ}(\vec R)$ to be the set of
$(\vec d_1,\ldots, \vec d_{\sQ})$ such that
$\vec d_1 <\cdots < \vec d_{ \sQ}$ and $\|\vec d_i\| \leq R_i$ for each $i$,
and we make the abbreviation $\vec\Delta_Q = (\vec d_1,\ldots, \vec d_Q)$.
We also use the abbreviations
\eqalign
\lbeq{T}
    T & = \prod_{q=1}^Q I[\{\vec D = \vec d_q \}\text{ on }
    \Vbold_{\vec d_q,\sM_{q-1}}^\smallsup{N}],
    \\
\lbeq{T3}
    T_3 & =
    \prod_{q=1}^{Q-1} I[\{\vec D = \vec d_q \}\text{ on }
    \Vbold_{\vec d_q,\sM_{q-1}}^\smallsup{N}]
    I[\{\vec D < \vec d_\sQ \}\text{ on } \Vbold_{\vec d_Q,\sM_{Q-1}}^\smallsup{N}]
    I[\vec D = \vec d_\sQ],
    \\
\lbeq{T4}
    T_4 & =
    \prod_{q=1}^Q I[\{\vec D = \vec d_q \}\text{ on }
    \Vbold_{\vec d_q,\sM_{q-1}}^\smallsup{N}]
      I[\|\vec D\| > R_{Q+1}],
    \\
\lbeq{T5}
    T_5 & = \prod_{q=1}^{\bar Q} I[\{\vec D = \vec d_q \}\text{ on }
    \Vbold_{\vec d_q,\sM_{q-1}}^\smallsup{N}]
    I[\vec D = \vec d_{\bar\sQ+1}].
\enalign
Iteration of \refeq{Iit} leads to
\eqalign
    & \sum_{\vec d_1 :\|\vec d_1\| \leq R_1}
     I[E^\smallsup{N} \text{ on } \Vbold_{\vec d_1,\sM_0}^\smallsup{N}]
    I[\vec D = \vec d_1]
    \nnb
    & \quad\quad =
     \sum_{Q=1}^{\bar Q}(-1)^{Q-1}
    \sum_{\vec \Delta_{ \sQ} \in \Scal_Q(\vec{R})}
    I[E^\smallsup{N} \text{ on } \Vbold_{\vec d_1,\sM_0}^\smallsup{N}]
    \big( T+T_3-T_4 \big)
    \nnb & \quad \quad \hspace{10mm}
\lbeq{Iitz}
    + (-1)^{\bar Q}
    \sum_{\vec \Delta_{\bar \sQ+1} \in \Scal_{\bar \sQ+1}(\vec{R})}
    I[E^\smallsup{N} \text{ on } \Vbold_{\vec d_1,\sM_0}^\smallsup{N}] T_5 .
\enalign

We insert the identity \refeq{Iitz}
into the right hand side of \refeq{PiMd}.
With \refeq{PiEE}, this gives
\eq
    \hat \Pi^\smallsup{N}
    =
    \Mcal^\smallsup{N}_{\vec \sR, \bar \sQ}
    +
    \Ecal^\smallsup{N}_{\vec \sR, \bar \sQ},
\en
where the right hand side is defined as follows.  First,
the main term is
\eqalign
\lbeq{Mcal}
    \Mcal^\smallsup{N}_{\vec \sR, \bar \sQ}
    &=
    \sum_{Q=1}^{\bar Q} (-1)^{Q-1}
    \sum_{\vec \Delta_{ \sQ} \in \Scal_Q(\vec{R})}
    \sum_{x,(u_j,v_j)}
    p^N  \Ebold^\smallsup{N}
    \Big[ I[E^\smallsup{N} \text{ on } \Vbold_{\vec d_1,\sM_0}^\smallsup{N}]
    T
     \Big].
\enalign
The error term is
\eq
    \Ecal^\smallsup{N}_{\vec \sR, \bar \sQ}
    =
    \Ecal^\smallsup{N}_{1,\sR_0}
    +
    \Ecal^\smallsup{N}_{2,\sR_0,\sR_1}
    +
    \sum_{j=3}^5 \Ecal^\smallsup{N}_{j,\vec \sR, \bar \sQ},
\en
where the first two terms on the right hand side are given by
\refeq{E1adef} and \refeq{E2def}, and
\eqalign
\lbeq{T3def}
    \Ecal^\smallsup{N}_{3, \vec \sR,\bar \sQ} & =
    \sum_{Q=1}^{\bar Q} (-1)^{Q-1}
    \sum_{\vec \Delta_{ \sQ} \in \Scal_Q(\vec{R})}
    \sum_{x,(u_j,v_j)}
    p^N  \Ebold^\smallsup{N}
    \Big[ I[E^\smallsup{N} \text{ on } \Vbold_{\vec d_1,\sM_0}^\smallsup{N}]
    T_3   \Big],
\\
\lbeq{T2def}
    \Ecal^\smallsup{N}_{4, \vec \sR,\bar \sQ} & =
    \sum_{Q=1}^{\bar Q}(-1)^{Q}
    \sum_{\vec \Delta_{ \sQ} \in \Scal_Q(\vec{R})}
    \sum_{x,(u_j,v_j)}
    p^N  \Ebold^\smallsup{N}
    \Big[ I[E^\smallsup{N} \text{ on } \Vbold_{\vec d_1,\sM_0}^\smallsup{N}]
    T_4  \Big],
\\
\lbeq{Udef}
    \Ecal^\smallsup{N}_{5, \vec \sR,\bar \sQ}
    &=
    (-1)^{\bar \sQ }
    \sum_{\vec \Delta_{\bar \sQ+1} \in \Scal_{\bar \sQ+1}(\vec R)}
    \sum_{x,(u_j,v_j)}
    p^N  \Ebold^\smallsup{N}
    \Big[ I[E^\smallsup{N} \text{ on } \Vbold_{\vec d_1,\sM_0}^\smallsup{N}]
    T_5 \Big].
\enalign

We will show in Section~\ref{sec-ee} that
$\bar{Q}(M)$ and $\vec R(M)$ can be chosen such that
$\Ecal^\smallsup{N}_{ \vec \sR,\bar \sQ} = O(\cn^{-M-1})$.
Given this bound on the error term,
to complete the proof of \refeq{PiNasy}
it suffices to prove the following
proposition.

\begin{prop}
Let $\bar{Q}(M)$ and  $\vec R(M)$ be given.
The main term $\Mcal^\smallsup{N}_{ \vec \sR, \bar \sQ}$
is a polynomial in $\cn$ and $p$ of the form
$\sum_{i,j=0}^{L_M} \alpha_{j,i,M}^\smallsup{N} \cn^j p^i$, as in
\refeq{PiNasy}, with rational coefficients
and with degree depending only on $M$.
\end{prop}

\proof
Since $\bar Q$ depends only on $M$, it suffices to show that each
term in the sum over $Q$ in \refeq{Mcal} is a polynomial of the desired
form.  Thus, for fixed $Q \in \{1,\ldots, \bar Q\}$, we will show that
\eq
\lbeq{PiMe}
    \sum_{\vec \Delta_{ \sQ} \in \Scal_Q(\vec{R})}
    \sum_{x,(u_j,v_j)}
    p^N  \Ebold^\smallsup{N}
    \Big[ I[E^\smallsup{N} \text{ on } \Vbold_{\vec d_1,\sM_0}^\smallsup{N}]
    \prod_{q=1}^Q I[\{\vec D = \vec d_q \}\text{ on }
    \Vbold_{\vec d_q,\sM_{q-1}}^\smallsup{N}]
     \Big]
\en
is a polynomial of the desired form.

We perform the sum over $\vec \Delta_\sQ$ in the order
\eq
    \sum_{\vec d_{ \sQ} : \| \vec d_\sQ\| \leq R_\sQ}
    \sumtwo{\vec \Delta_{ \sQ-1} \in \Scal_{\sQ-1}(\vec{R}):}
    {\vec d_{\sQ-1} < \vec d_\sQ },
\en
which puts the sum over $\vec d_\sQ$ last.
When $Q=1$, the second sum is absent.
We define an equivalence relation on $\{0,1\}^n$ by regarding
$\vec d $ and $\vec d'$ as equivalent
if $\| \vec d\| = \|\vec d'\|$,
i.e., if $\vec d$ and $\vec d'$ have the same number of components
taking the value 1.  If $\vec d$ and $\vec d'$ are equivalent,
then by symmetry
they give rise to equal contributions to the sum over $\vec d_\sQ$
in \refeq{PiMe}.  Let $[\vec d]_m$ denote the element of
$\{0,1\}^n$ that consists of $m$ ones followed by $n-m$ zeros.
Then \refeq{PiMe} is equal to
\eq
\lbeq{polya}
    \sum_{m=1}^{R_\sQ} {n \choose m}
    \sumtwo{\vec \Delta_{ \sQ-1} \in \Scal_{\sQ-1}(\vec{R}):}
    {\vec d_{\sQ-1} < [\vec d]_m }
    \sum_{x,(u_j,v_j)}
    p^N  \Ebold^\smallsup{N}
    \Big[ I[E^\smallsup{N} \text{ on } \Vbold_{\vec d_1,\sM_0}^\smallsup{N}]
    \prod_{q=1}^Q I[\{\vec D = \vec d_q \}\text{ on }
    \Vbold_{\vec d_q,\sM_{q-1}}^\smallsup{N}]
     \Big],
\en
where $\vec d_\sQ$ is equal to $[\vec d]_m$.

The number of terms in the sum over $\vec \Delta_{\sQ-1}$ depends only
on $\vec R$, and hence depends only on $M$.  Also,
the cardinality of the set
$\Vbold_{[\vec d]_m,\sM_{Q-1}}^\smallsup{N}$ is bounded by a constant
depending only on $M$, uniformly in $m \leq R_\sQ$, and this set
contains the sets $\Vbold_{\vec d_q,\sM_{q-1}}^\smallsup{N}$
for $q < Q$.
The event
$\{E^\smallsup{N} \text{ on } \Vbold_{\vec d_1,\sM_0}^\smallsup{N}\}$
implies that
$x, u_j, v_j\in \Vbold_{\vec d_1,\sM_0}^\smallsup{N}$,
and hence the number of terms in the sum over $x,(u_j,v_j)$
is bounded by a constant depending only on $M$.
Thus, it suffices to prove that the expectation in \refeq{polya} is
a polynomial of the desired form.

The expectation in \refeq{polya}
is the probability of an event that only depends on the
occupation status of bonds in $\Vbold_{[\vec d]_m,\sM_{\sQ-1}}^\smallsup{N}$.
Explicitly, this probability is the sum, over the
finitely many configurations on $\Vbold_{[\vec d]_m,\sM_{\sQ-1}}^\smallsup{N}$
for which the product of indicators is $1$, of $p^x(1-p)^y$,
where $x$ and $y$ are the number of occupied and vacant bonds,
respectively, in the configuration.  Thus, this probability is a polynomial
in $p$ with integer coefficients.
Therefore, taking into account the binomial coefficient
in \refeq{polya}, the coefficients of
the polynomial \refeq{polya}
in $n$ and $p$ must all be rational numbers.
(We note that rational coefficients arise here rather than the integer
coefficients found for the self-avoiding walk in \cite{HS95}, because in
counting dimensions under symmetry, combinations arise here
whereas for the self-avoiding walk it was permutations.)
\qed

\subsection{Error estimates}
\label{sec-ee}

We now formulate three lemmas which show that
we can choose $\bar{Q}$ and $R_0<R_1<\cdots<R_{\bar \sQ+1}$,
all depending on $M$,
such that each of
$\Ecal_{1,\sR}^\smallsup{N}$, $\Ecal_{2,\sR,\sR_1}^\smallsup{N}$,
and $\Ecal^\smallsup{N}_{j,\vec \sR,\bar \sQ}$ ($j=3,4,5$)
is $O(\cn^{-M-1})$.  As noted above, these estimates imply
\refeq{PiNasy} and thus complete the proof of Proposition~\ref{prop-Piae}.
Proofs of the three lemmas will be given in
Section~\ref{sec-error}.

The first lemma gives the desired bound on $\Ecal_{1,\sM_0}^\smallsup{N}$.

\begin{lemma}
\label{lem-E1}
Let $M \geq 1$.  There exists a $K_1=K_1(M)$ and an $r_0(M)$ such that
if $R_0 \geq r_0(M)$ then for $N \leq M$ and $p \leq \bar p_c$,
\eq
    |\Ecal_{1,\sM_0}^\smallsup{N}|
    \leq K_1 \cn^{-M-1}.
\en
\end{lemma}

The corollary to the
second lemma gives the desired bound on
$\Ecal^\smallsup{N}_{3,\sR, \vec \sR,\bar \sQ}$.

\begin{lemma}
\label{lem-E3}
Let $M \geq 1$.  There exists a $K_2=K_2(M)$ and an $r(M)$ such that
if $R \geq r(M)$ then for $N \leq M$ and $p \leq \bar p_c$,
\eqalign
\lbeq{E3lem}
    &
    \sum_{\vec d \in \{0,1\}^n}
    \sum_{x,(u_j,v_j)}
    p^N  \Ebold^\smallsup{N}
    \left[ I[F^\smallsup{N}]
    I[\{\vec D < \vec d\}
    \text{ on } \Vbold_{\vec d,\sM}^\smallsup{N}]
    I[\vec D = \vec d] \right]
    \leq K_2 \cn^{-M-1}.
\enalign
\end{lemma}

\begin{cor}
\label{cor-T3}
Let $M\geq 1$ and $\bar Q \geq 2$.  If
$r(M) \leq R_{\sQ} \leq R_{\bar\sQ}$ for $Q=0,\ldots,\bar Q-1$,
then for $N \leq M$ and $p \leq \bar p_c$,
    \eq
    |\Ecal^\smallsup{N}_{3, \vec \sR,\bar \sQ}|
    \leq \bar Q 2^{(\bar Q-1)R_{\bar \sQ}} K_2 \cn^{-M-1}.
    \en
\end{cor}

\proof
By \refeq{Fbdon} and \refeq{T3},
the summand in the sum over $Q$ in \refeq{T3def} is bounded above by
    \eq
    \lbeq{cT3a}
    \sum_{\vec \Delta_{\sQ} \in \Scal_{ \sQ}(\vec R)}
    \sum_{x,(u_j,v_j)}
    p^N  \Ebold^\smallsup{N} \left[ I[F^\smallsup{N} ]
    I[\{\vec D < \vec d_\sQ \}\text{ on }
    \Vbold_{\vec d_Q,\sM_{Q-1}}^\smallsup{N}]
    I[\vec D = \vec d_\sQ]\right].
    \en
Given $\vec d_{\sQ}$ with $\| \vec d_\sQ\| \leq R_\sQ \leq R_{\bar\sQ}$,
the number of $\vec d_1,\ldots,  \vec d_{Q-1}$
with $\vec d_1 <\cdots <  \vec d_{Q-1} < \vec d_{ Q}$ and
$\|\vec d_i\| \leq R_i$ is bounded above by
$2^{(\bar Q-1)R_{\bar \sQ}}$.
Thus \refeq{cT3a} is bounded above by
\eq
    2^{(\bar Q-1)R_{\bar \sQ}}
    \sum_{\vec d_\sQ \in \{0,1\}^n}
    \sum_{x,(u_j,v_j)}
    p^N  \Ebold^\smallsup{N} \left[ I[F^\smallsup{N} ]
    I[\{\vec D < \vec d_\sQ \}\text{ on }
    \Vbold_{\vec d_Q,\sM_{Q-1}}^\smallsup{N}]
    I[\vec D = \vec d_\sQ]\right],
\en
and now \refeq{E3lem} can be applied.
\qed

The corollary to the
third lemma gives the desired bounds on
$\Ecal_{2,\sR,\sR_1}^\smallsup{N}$,
and $\Ecal^\smallsup{N}_{j, \vec \sR,\bar\sQ}$ ($j=4,5$).

\begin{lemma}
\label{lem-E2}
Let $M \geq 1$ and $R \geq 1$.
There exists a constant $C(R,M)$, and
a sequence $g_R$ with $\lim_{R\to\infty} g_R = \infty$,
such that for $N \leq M$ and $p \leq \bar p_c$,
\eqalign
\lbeq{E2lema}
    \sum_{x,(u_j,v_j)}
    p^N  \Ebold^\smallsup{N}
    \left[ I[F^\smallsup{N} ]
    I[\|\vec D\| > R] \right]
    &
    \leq
    C(R,M) \cn^{-g_R}.
\enalign
\end{lemma}

\begin{cor}
\label{cor-T2}
For $M\geq 1$, we can choose $\bar{Q}=\bar{Q}(M)$ and
$R_i=R_i(M)$ ($i=1,\ldots,\bar Q +1$), with
$R_1<R_2<\cdots <R_{\bar\sQ+1}$,
such that for $N \leq M$, $p \leq \bar p_c$, and for any $R_0$,
    \eq
    \left|\Ecal_{2,\sR_0,\sR_1}^\smallsup{N}\right| \leq C(R_1,M) \cn^{-M-1},
    \quad
    \left| \Ecal^\smallsup{N}_{4, \vec \sR,\bar\sQ}\right|
    \leq C'(\bar Q, \vec R)\cn^{-M-1},
    \en
    \eq
    \left| \Ecal^\smallsup{N}_{5, \vec \sR,\bar\sQ}\right|
    \leq
    2^{\bar QR_{\bar\sQ +1}}
    C(R_{\bar Q},M) \cn^{-M-1},
    \en
where $C'(\bar Q, \vec R)$ is a constant depending on
$\bar Q$ and $\vec R$, and hence only
on $M$.
\end{cor}

\proof
The bound on $\Ecal_{2,\sR_0,\sR_1}^\smallsup{N}$ is an immediate consequence
of \refeq{E2def}, \refeq{Fbdon} and \refeq{E2lema}.

For the bound on $\Ecal^\smallsup{N}_{5, \vec \sR,\bar\sQ}$, we note
from \refeq{Udef}, \refeq{Fbdon} and \refeq{T5} that
\eq
     \left| \Ecal^\smallsup{N}_{5, \vec \sR,\bar\sQ} \right|
    \leq
    \sum_{\vec \Delta_{\bar \sQ+1} \in \Scal_{\bar \sQ+1}(\vec R)}
    \sum_{x,(u_j,v_j)}
    p^N  \Ebold^\smallsup{N} \left[ I[F^\smallsup{N} ]
    I[\vec D = \vec d_{\bar Q+1}]\right].
\en
For fixed $\vec d_{\bar\sQ+1}$, the
number of $\vec d_1,\ldots,  \vec d_{\bar Q}$
with $\vec d_1 <\cdots <  \vec d_{\bar Q} < \vec d_{\bar Q +1}$ and
$\|\vec d_i\| \leq R_i \leq R_{\bar \sQ+1}$
is bounded above by $2^{{\bar Q} R_{\bar Q +1}}$.
On the other hand, since the $\vec d_i$ are strictly
increasing, it must be the
case that $\|\vec d_{\bar Q+1}\| \geq \bar Q $, and hence
\eq
    \left| \Ecal^\smallsup{N}_{5, \vec \sR,\bar\sQ} \right|
    \leq
    2^{{\bar Q} R_{\bar Q +1}}
    \sum_{x,(u_j,v_j)}
    p^N  \Ebold^\smallsup{N} \left[ I[F^\smallsup{N} ]
    I[\|\vec D\| \geq \bar Q ]\right].
\en
By Lemma~\ref{lem-E2},
we can choose $\bar Q =\bar Q(M)$
such that
\eq
    \left|\Ecal^\smallsup{N}_{5, \vec \sR,\bar\sQ} \right|
    \leq
    2^{{\bar Q} R_{\bar Q +1}}
    C(R_{\bar Q},M) \cn^{-M-1}.
\en

Finally, we prove the bound on $\Ecal^\smallsup{N}_{4, \vec \sR,\bar\sQ}$.
The number of $\vec d_\sQ$ with $\|\vec d_{\sQ}\| \leq R_\sQ$ is
at most $n^{\|\vec d_Q\|} \leq n^{R_Q}$.
Given such a $\vec d_{\sQ}$, the number of
$(\vec d_1,\ldots, \vec d_{\sQ-1})$
with $\vec d_1 < \cdots < \vec d_{\sQ -1} < \vec d_{\sQ}$ is
at most $2^{(Q-1)R_{\sQ}}$.
Therefore,
\eq
\lbeq{T2bd}
    \left| \Ecal^\smallsup{N}_{4,\vec \sR,\bar\sQ} \right|
    \leq
    \sum_{Q=1}^{\bar Q}
    2^{(Q-1)R_{\sQ}} n^{R_Q}
    \sum_{x,(u_j,v_j)}
    p^N  \Ebold^\smallsup{N}
    \left[ I[F^\smallsup{N} ]
    I[\|\vec D\| > R_{Q+1}] \right].
\en
By Lemma~\ref{lem-E2}, given any choice of $R_0$ and any choice of $R_1$,
we can choose $R_{\sQ+1}=R_{\sQ+1}(R_\sQ,M)>R_\sQ$ sequentially and increasing
for $Q=1,\ldots, \bar Q$, so that
\eqalign
    \left| \Ecal^\smallsup{N}_{4,\vec \sR,\bar\sQ} \right|
    & \leq
    \sum_{Q=1}^{\bar Q}
    2^{(Q-1)R_{\sQ}} n^{R_Q}
    C(R_\sQ,M)\cn^{-M-1-R_Q}
    \nnb &
    \leq
    \cn^{-M-1}\sum_{Q=1}^{\bar Q}
    2^{(Q-1)R_{\sQ}}
    C(R_\sQ,M)\cn^{-M-1}.
\enalign
This is the desired estimate.  It is at this point that we make use of the
flexibility to choose a sequence $R_1,R_2,\ldots$.
\qed

\vskip 0.5cm

This gives the desired bounds on the error terms.
The value of $\bar Q(M)$ is fixed by Corollary~\ref{cor-T2},
we take $R_0=r_0(M)\vee r(M)$, and
we fix $R_{\bar\sQ+1}> \cdots >R_1>R_0$ according to Corollary~\ref{cor-T2}.
Then the restrictions of Lemma~\ref{lem-E1}
and Corollaries~\ref{cor-T3} and \ref{cor-T2} are all obeyed.

It remains to
prove Lemmas~\ref{lem-E1}, \ref{lem-E3} and \ref{lem-E2}.
This will be done in Section~\ref{sec-error}.

\section{Proof of error estimates}
\label{sec-error}

In this section, we complete the proof of Proposition~\ref{prop-Piae}
by proving Lemmas~\ref{lem-E1}, \ref{lem-E3} and \ref{lem-E2}.
We begin by recalling some basic facts.

\subsection{Preliminaries}
\label{sec-pre}

Let $D(y-x) = \cn^{-1}$ if $x$ and $y$ are neighbours,
and $D(y-x)=0$ otherwise.  Thus $D(y-x)$ is the transition
probability for simple random walk on $\gr$ to make a step
from $x$ to $y$.  Let
$\tau_p(y-x) = \Pbold_p(x \conn y)$ denote the two-point
function, and let
$\tau_p^\smallsup{i}(x)$ denote the probability that there is
an occupied (self-avoiding) path from $0$ to $x$ of length at least $i$.

We define the Fourier transform of
an absolutely summable function $f$ on the vertex set $\Vbold$
of $\gr$ by
\eq
    \hat{f}(k) = \sum_{x \in \Vbold} f(x) e^{ik\cdot x}
    \quad \quad
    (k \in \Vbold^*),
\en
where $\Vbold^* = \{0,\pi\}^n$ for $\qn$ and $\Vbold^* = [-\pi,\pi]^n$
for $\Z^n$.
Let
\eq
    (f*g)(x) = \sum_{y \in \Vbold} f(y)g(x-y)
\en
denote convolution.
Recall from \cite{AN84} that $\hat\tau_p(k) \geq 0$ for all $k$.

For $i,j$ non-negative integers, let
\eqalign
    T_p^\smallsup{i,j}
    & =
     \begin{cases}
    2^{-n} \sum_{k \in \{0,\pi\}^n} |\hat{D}(k)|^i \hat{\tau}_p(k)^j
    & (\gr = \qn)
    \\
    \int_{[-\pi,\pi]^n} |\hat{D}(k)|^i \hat{\tau}_p(k)^j \frac{d^nk}{(2\pi)^n}
    & (\gr = \Z^n),
    \end{cases}
    \\
    T_p &= \sup_{x} (p\cn)(D*\tau_p*\tau_p*\tau_p)(x).
\enalign
Recall from \cite[Section~3]{HS04a} that
for $\gr = \Z^n$ and $\gr = \qn$,
there are constants $K_{i,j}$ and $K$
such that
for all $p \leq \bar p_c(\gr)$,
\eqalign
\lbeq{Tij-bd}
     T_p^\smallsup{i,j}
     & \leq K_{i,j} \cn^{-i/2} \quad \text{($i,j\geq 0$)},
     \\
\lbeq{Tpbound}
     T_p & \leq K \cn^{-1},
\\
\lbeq{tauiT}
    \sup_x \tau_p^\smallsup{i}(x)
    & \leq
    \begin{cases}
    K \cn^{-1} & (i =1) \\
    2^i K_{i,1}\cn^{-i/2} & (i \geq 2).
    \end{cases}
\enalign
The above bounds are valid for $n \geq 1$ for $\qn$, and for $n$ larger
than an absolute constant for $\Z^n$, except \refeq{Tij-bd} also requires
$n \geq 2j+1$ for $\Z^n$.

%
%
%
%

\subsection{Proof of Lemmas~\ref{lem-E1}, \ref{lem-E3} and \ref{lem-E2}}
\label{sec-error2}
\newcommand{\sMR}{{\sss M}}

Now we prove
Lemmas~\ref{lem-E1}, \ref{lem-E3} and \ref{lem-E2}.
The proofs use the following proposition.
Recall the definition of $\Pcal_j$ in Definition~\ref{def-bb}(ii),
and let ${\cal P}_{j,{\sss L}}$ denote the
subset of ${\cal P}_j$ consisting of paths of length less than $L$.

\begin{prop}
\label{prop-PibdsM}
Let $M \geq 1$.
There is a constant $K_3=K_3(M)$ such that for   $N \leq M$
and $p \leq \bar p_c$,
    \eq
    \lbeq{Pibd1}
    \sum_{x,(u_j,v_j)}
    p^N  \Ebold^\smallsup{N} \left[ I[F^\smallsup{N} ]
    I[\{\exists \text{ occupied }\omega \in {\cal P}_{j}\backslash
    \Pcal_{j,10(M+1)} \}]\right] \leq K_3\cn^{-M-1} .
    \en
\end{prop}

\proof
We assume some familiarity with diagrammatic estimates,
as in \cite[Section~4]{BCHSS04b}.

We begin by rewriting \refeq{FN-def} as
    \eq
    \lbeq{FFj}
    F^\smallsup{N}=\bigcup_{\vec{t},\vec{w}, \vec{z}}\bigcap_{j=0}^{N} F_j,
    \en
where we have made the abbreviations
$F_0=F_0(0, u_0, w_0, z_0)_0$,
$F_i=F(v_{i-1},t_i,z_i,u_i, w_i, z_{i+1})_i$, and
$F_N=F_N(v_{N-1},t_N,z_N,x)_N$.
We will use the estimate
    \eq
    \lbeq{PGF}
    \Pbold^\smallsup{N}\Big(G\cap \bigcup_{\vec{t},\vec{w}, \vec{z}}
    \bigcap _{j=0}^{N} F_j\Big)
    \leq \sum_{\vec{t},\vec{w}, \vec{z}} \Pbold^\smallsup{N}
    \Big(G\cap \bigcap _{j=0}^{N} F_j\Big),
    \en
with $G = \{\exists \text{ occupied }\omega \in {\cal P}_{j}\backslash
\Pcal_{j,10(M+1)} \}$ (for fixed $j,M$).
The standard bounds on $\hat\Pi^\smallsup{N}_p$
use \refeq{PGF} with
$G$ equal to the whole probability space.  In the standard bounds,
after applying the BK inequality, each of the disjoint
connections in \refeq{Fdefa}--\refeq{Fdefsa}, say $\{y_1\conn y_2\}$,
gives rise to a two-point function $\tau_p(y_2-y_1)$.
The overall effect is to bound $\hat\Pi^\smallsup{N}_p$ by a Feynman diagram,
which is then bounded by products of $T^\smallsup{i,3}$
as in \cite[Section~4]{BCHSS04b}.  We will
modify this standard procedure to prove
the proposition.

We decompose $G$ as $G=G_{1}\cup G_{2}$, where $G_{1}$
is the event that one of the disjoint connections in
\refeq{Fdefa}--\refeq{Fdefsa}, say $\{y_1\conn y_2\}$, is replaced by
the disjoint
occurrence of $\{y_1\conn y_2$ via a path consisting
of at least $2(M+1)$ occupied bonds$\}$,
and $G_{2}=G\backslash G_{1}.$

For the contribution due to $G_{1}$, the standard diagrammatic
bounds give an upper bound identical to that for $\hat \Pi^{\smallsup{N}}_p$,
except that the factor $\tau_p(y_2-y_1)$ is replaced by
$\tau_p^\smallsup{2M+2}(y_2-y_1)$.
This replaces the bound
    \eq
    \hat\Pi^\smallsup{N}_p
    \leq
    \begin{cases}
    T_p  & (N=0)
    \\
    T_p^\smallsup{0,3}(2 T_p^\smallsup{0,3}T_p)^{N } & (N\geq 1)
    \end{cases}
    \lbeq{PiNbd}
    \en
of \cite[Proposition~4.1]{BCHSS04b}
by a sum of terms in which one factor
 $T^\smallsup{0,3}_p$ or $T_p$
is replaced by $T^{\smallsup{2M+2,3}}_p$.
The number of terms in the sum depends only on $M$.
By \refeq{Tij-bd}, this new upper bound
is at most
$O(\cn^{-M-1})$, as required. Thus, we are left to deal with $G_{2}$.

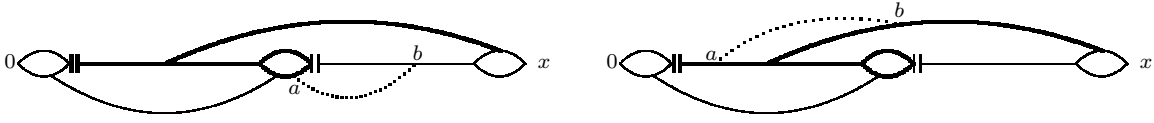
\begin{figure}[t]
\vskip2cm
\begin{center}
\setlength{\unitlength}{0.00450in}%
\begin{picture}(1400,-100)
\put(-5,15){${\scriptstyle 0}$}
\put(615,15){${\scriptstyle x}$}
\thinlines
\qbezier(10,20)(40,-10)(70,20)
\qbezier(10,20)(40,50)(70,20)
\qbezier(50,5)(180,-80)(310,5)
\Thicklines
\put(73,10){\line(0,1){20}}
\put(80,10){\line(0,1){20}}
\put(80,20){\line(1,0){210}}
\qbezier(290,20)(320,-10)(350,20)
\qbezier(290,20)(320,50)(350,20)
\qbezier(180,20)(375,110)(570,35)
\thicklines
\qbezier[25](330,5)(400,-50)(470,18)
\put(325,-15){${\scriptstyle a}$}
\put(470,25){${\scriptstyle b}$}
\thinlines
\put(352,10){\line(0,1){20}}
\put(360,10){\line(0,1){20}}
\put(360,20){\line(1,0){180}}
\qbezier(540,20)(570,-10)(600,20)
\qbezier(540,20)(570,50)(600,20)
\put(695,15){${\scriptstyle 0}$}
\put(1315,15){${\scriptstyle x}$}
\thinlines
\qbezier(710,20)(740,-10)(770,20)
\qbezier(710,20)(740,50)(770,20)
\qbezier(750,5)(880,-80)(1010,5)
\Thicklines
\put(773,10){\line(0,1){20}}
\put(780,10){\line(0,1){20}}
\put(780,20){\line(1,0){210}}
\qbezier(990,20)(1020,-10)(1050,20)
\qbezier(990,20)(1020,50)(1050,20)
\qbezier(880,20)(1075,110)(1270,35)
\thicklines
\qbezier[25](820,20)(900,90)(1030,67)
\put(810,25){${\scriptstyle a}$}
\put(1030,75){${\scriptstyle b}$}

\thinlines
\put(1052,10){\line(0,1){20}}
\put(1060,10){\line(0,1){20}}
\put(1060,20){\line(1,0){180}}
\qbezier(1240,20)(1270,-10)(1300,20)
\qbezier(1240,20)(1270,50)(1300,20)
\end{picture}
\end{center}
\vskip1cm

\caption{Examples of disjoint connections satisfying $G_2$.
Each solid path joining a pair of vertices has length at most $2(M+1)$,
and the first and second dotted paths have lengths at least
$2(M+1)$ and $6(M+1)$, respectively.}
\label{fig-i}
\end{figure}

To estimate the right hand side of \refeq{PGF} with $G=G_2$,
we may assume that each of the connections in
$F_0,\ldots, F_N$ is achieved by a path consisting
of at most $2(M+1)$ bonds.  On the other hand, there must exist
an occupied path in $\Pcal_j$ consisting of at least $10(M+1)$
bonds.  The coexistence of this path with the disjoint occupied
paths required by the event $F_j$ implies that we can find
vertices $a,b$ and
disjoint paths such that two of the disjoint connections
in $F_j$ and/or $F_{j+1}$, say $\{y_1\conn y_2\}$
and $\{y_3 \conn y_4\}$,  become replaced by
\eq
\lbeq{yab}
    \{y_1\conn a\} \circ \{a\conn y_2\} \circ \{a\conn b\} \circ
    \{y_3\conn b\} \circ \{b\conn y_4\} .
\en
Moreover, the connection from $a$ to $b$ must be achieved by a path
of length at least $2(M+1)$.  This follows from the fact that
the long occupied path in $\Pcal_j$ can partially
coincide with at most four of the paths realizing the
connections of $F_j$, and these paths have total length
at most $8(M+1)$.  See Figure~\ref{fig-i}.  Therefore, we can
in fact replace \refeq{yab} by
\eq
\lbeq{yablong}
    \{y_1\conn a\} \circ \{a\conn y_2\} \circ \{a\conn b
    \text{ by a path of length at least $2(M+1)$}\} \circ
    \{y_3\conn b\} \circ \{b\conn y_4\} .
\en

In an upper bound achieved via the BK inequality, the two
factors  $\tau_p(y_1-y_2)\tau_p(y_3-y_4)$ normally present in the upper
bound on $\hat{\Pi}^\smallsup{N}_p$
are replaced by
    \eqalign
    &
    \sum_{a,b} \tau_p(y_1-a)\tau_p(a-y_2)
    \tau^\smallsup{2M+2}_p(b-a)\tau_p(y_3-b)\tau_p(b-y_4)
    \nnb
    & \qquad \leq
    \left(\sup_{a,b} \tau^\smallsup{2M+2}_p(b-a)\right)
    \sum_{a} \tau_p(y_1-a)\tau_p(a-y_2)\sum_{b} \tau_p(y_3-b)\tau_p(b-y_4)
    \nnb
    &\qquad \leq 2^{2M+2}K_{2M+2,1}\Omega^{-M-1}
    (\tau_p*\tau_p)(y_1-y_2)(\tau_p*\tau_p)(y_3-y_4)
    ,
    \enalign
using \refeq{tauiT} to estimate
$\sup_{a,b} \tau^\smallsup{2M+2}_p(b-a)$.
It remains to show that the result of replacing the factors
$\tau_p(y_1-y_2)\tau_p(y_3-y_4)$
by
$(\tau_p*\tau_p)(y_1-y_2)(\tau_p*\tau_p)(y_3-y_4)$,
in the standard bound on $\hat\Pi^\smallsup{N}_p$,
gives rise to a bounded quantity.

This replacement has the effect of adding a vertex to each of
the lines joining $y_1$ to $y_2$ and $y_3$ to $y_4$ in the standard
diagrammatic bound on $\hat\Pi^\smallsup{N}_p$.
This leads to a bound in which one factor of
$T^{\smallsup{0,3}}_p$ or $T_p$ in \refeq{PiNbd}
is replaced by $T^{\smallsup{0,5}}_p$ or $T^{\smallsup{1,5}}_p$, or
a product of two
factors of $T^{\smallsup{0,3}}_p$ and/or $T_p$
is replaced by $T^{\smallsup{0,4}}_p
T^{\smallsup{0,4}}_p$ (taking an upper bound),
depending on whether the two vertices
are added to the same triangle or not.
The result is finite by \refeq{Tij-bd}.
\qed

\vskip 0.5 cm
\noindent
{\em Proof of Lemma \ref{lem-E1}.}
We write simply $R$ in place of $R_0$, and define
    \eqalign
    \tilde C_{j,\sM}
    &=
    \{y:\{ v_{j-1} \conn y
    \text{ without using $(u_j,v_j)\}$ on }\Bbold_\sM^\smallsup{N}\}
    \quad (j=0,\ldots, N-1),
    \\
    E_{j,\sR}& =
    \begin{cases}
    E_0 & (j=0) \\
    E'(v_{j-1}, u_j; \tilde C_{j-1,\sM})_j & ( j = 1,\ldots, N),
    \end{cases}
    \\
    E^\smallsup{N}_\sR
    &= \bigcap_{j=0}^N  E_{j,\sR}.
    \enalign
By definition,
    \eq
    I[E^\smallsup{N} \text{ on } \Bbold_\sM^\smallsup{N}]
    =I[E^\smallsup{N}_\sR \text{ on } \Bbold_\sM^\smallsup{N}],
    \en
and we may therefore rewrite the difference occurring in the definition
of $\Ecal_{1,\sR}^\smallsup{N}$ in \refeq{E1adef} as
    \eq
    \lbeq{basrewon}
    I[E^\smallsup{N}]
    -
    I[E^\smallsup{N} \text{ on } \Bbold_\sM^\smallsup{N}]
    =
    \big(
    I[E^\smallsup{N}_\sR]
    -
    I[E^\smallsup{N}_\sR \text{ on } \Bbold_\sM^\smallsup{N}]
    \big)
    +
    \big(
    I[E^\smallsup{N}]
    -
    I[E^\smallsup{N}_\sR]
    \big).
    \en

We will prove that
\eqalign
\lbeq{1stdif}
    \big|I[E^\smallsup{N}_\sR]
    -
    I[E^\smallsup{N}_\sR \text{ on } \Bbold_\sM^\smallsup{N}]\big|
    & \leq
    2I[F^\smallsup{N}] \sum_{j=0}^N
    I[ 
    B_j \cap (\Bbold_\sM^\smallsup{N})^c
    \neq \varnothing],
    \\
\lbeq{2nddif}
    \big|I[E^\smallsup{N}]
    -
    I[E^\smallsup{N}_\sR]
    \big|
    & \leq
    I[F^\smallsup{N}]
    \sum_{j=1}^N
    I[
    B_j\cap \tilde C_{j-1}\backslash \tilde C_{j-1,\sR}\neq \varnothing].
\enalign
Assuming \refeq{1stdif}--\refeq{2nddif},
the proof is completed as follows.
Consider first a configuration contributing to the
summand on the right hand side of \refeq{1stdif}.  In such a configuration,
there must be a backbone path of length
$R/(N+1) \geq R/(M+1)$, since otherwise
no backbone path could exit $\Bbold_\sM^\smallsup{N}$.
We take $R = 10(M+1)^2$ and apply Proposition~\ref{prop-PibdsM}.
Similarly, in a configuration contributing to the summand in the
right hand side of \refeq{2nddif}, if there is an element of
$B_j \cap \tilde C_{j-1}\backslash \tilde C_{j-1,\sR}$ that
lies outside of $\Bbold_\sM^\smallsup{N}$, then
there must be a path in $\Pcal_j$ that exits $\Bbold_\sM^\smallsup{N}$.
This implies that there is an occupied path in some $\Pcal_i$ of length
at least $R/(M+1)$, and we again take $R = 10(M+1)^2$ and apply
Proposition~\ref{prop-PibdsM}.
On the other hand,
if $B_j \cap \tilde C_{j-1}\backslash \tilde C_{j-1,\sR}
\subset \Bbold_\sM^\smallsup{N}$, then a path in
$\tilde{C}_{j-1}$ must travel from $v_{j-2}$ outside of
$\Bbold_\sM^\smallsup{N}$ before intersecting $B_j$, so there
must then be a path in $\Pcal_{j-1}$ that exits $\Bbold_\sM^\smallsup{N}$,
and again the previous argument applies.
It remains to prove \refeq{1stdif}--\refeq{2nddif}.

\smallskip \noindent \emph{Proof of \refeq{1stdif}.}
By \refeq{on/inprop1},
\eqalign
\lbeq{Iona}
    I[E^\smallsup{N}_\sR]
    -
    I[E^\smallsup{N}_\sR \text{ on } \Bbold_\sM^{\smallsup{N}}]
    & =
    I[E^\smallsup{N}_\sR \cap \{(E^\smallsup{N}_\sR)^c \text{ on }
    \Bbold_\sM^\smallsup{N}\}]
     -
    I[(E^\smallsup{N}_\sR)^c \cap \{E^\smallsup{N}_\sR\text{ on }
    \Bbold_\sM^\smallsup{N}\}].
\enalign
By \refeq{Fbdon}, $\{E^\smallsup{N}\text{ on } \Bbold_\sM^\smallsup{N}\} \subset
F^\smallsup{N}$, and the proof of \refeq{Fbda} also easily extends to yield
$E^\smallsup{N}_\sR \subset F^\smallsup{N}$.
By \refeq{on/inprop1}--\refeq{on/inprop2},
it follows that
    \eqalign
    I[E^\smallsup{N}_\sR \cap \{(E^\smallsup{N}_\sR)^c\text{ on }
    \Bbold_\sM^\smallsup{N}\}]
    & \leq I[E^\smallsup{N}_\sR] \sum_{j=0}^N I[E_{j,\sR}\cap
    \{E_{j,\sR}^c \text{ on } \Bbold_\sM^\smallsup{N}\}]
    \nnb & \leq
    I[F^\smallsup{N}]
    \sum_{j=0}^N I[E_{j,\sR}\cap
    \{E_{j,\sR}^c \text{ on } \Bbold_\sM^\smallsup{N}\}],
    \\
    I[(E_\sR^\smallsup{N})^c\cap
    \{E_\sR^\smallsup{N}\text{ on } \Bbold_\sM^\smallsup{N}\}]
    &\leq I[\{E_\sR^\smallsup{N}\text{ on } \Bbold_\sM^\smallsup{N}\}]
    \sum_{j=0}^N I[E_{j,\sR}^c\cap \{E_{j,\sR}\text{ on }
    \Bbold_\sM^{\smallsup{N}}\}]
    \nnb & \leq
    I[F^\smallsup{N}]
    \sum_{j=0}^N I[E_{\sR,j}^c\cap \{E_{j,\sR}\text{ on }
    \Bbold_\sM^{\smallsup{N}}\}].
    \enalign
Thus, it suffices to show that (i) if $E_{j,\sR}\cap
    \{E_{j,\sR}^c$ on $\Bbold_\sM^\smallsup{N}\}$ occurs,
or (ii) if $E_{j,\sR}^c\cap \{E_{j,\sR}$ on
$\Bbold_\sM^{\smallsup{N}}\}$ occurs,
then there is a path in $B_j$ that
exits $\Bbold_\sM^{\smallsup{N}}$.

We first consider case (i).  It is clear that if $E_0 \cap
\{E_0^c \text{ on } \Bbold_\sM^{\smallsup{N}}\}$
occurs, then there is a path in $B_0$
that exits $\Bbold_\sM^{\smallsup{N}}$.
So we consider $j \geq 1$.  Given sets $A,B$ of vertices,
it suffices to show that if $E'(v,x;A)\cap \{E'(v,x;A)^c$ on $B\}$ occurs,
then there must be an occupied path from $v$ to $x$ that exits $B$.
Recall from \refeq{317} that the event $E'(v,x;A)$ is the intersection
of the event $\{v \ct{A} x\}$ with the NP condition.
Therefore,
    \eqalign
    \!\! I[E'(v,x;A)\cap \{E'(v,x;A)^c \text{ on } B\}]
    &=I[E'(v,x;A)]I[\{v \ct{A} x\}^c \text{ on $B$}]\nnb
    &\quad +I[E'(v,x;A)]I[\{v \ct{A} x\} \text{ on $B$}]
    I[{\rm NP}^c \text{ on $B$}].
    \enalign
In the first term on the right hand side, the event $E'(v,x;A)$
requires that $\{v \ct{A} x\}$.
If every path from $v$ to $x$ stays inside $B$, then also
$\{v \ct{A} x\}$ on $B$.
The second factor therefore ensures
that there must be a connection from $v$ to $x$ that exits $B$,
as required.  In the second term on the right hand side, NP holds,
but not on $B$.  This means that, on $B$, there is a pivotal bond $(u',v')$
for the connection from $v$ to $x$ such that $v \ct{A} u'$,
but that there is no such bond when the entire
configuration on $\gr$ is used.
This can only happen if there is an occupied path from $v$ to $x$
that exits $B$, with this path either making $(u',v')$ no longer pivotal,
or providing a path from $v$ to $u'$ that does not intersect $A$.

Next, we consider case (ii).
The case $j=0$ cannot occur.  We consider $j \geq 1$, and proceed
as in the proof of case (i).
By \refeq{317},
    \eqalign
    I[E'(v,x;A)^c\cap \{E'(v,x;A) \text{ on }B\}]
    &
    =I[E'(v,x;A) \text{ on } B]
    I[\{v \ct{A} x\}^c \cap {\rm NP}]
    \nnb &\quad
    +I[E'(v,x;A) \text{ on } B]
    I[{\rm NP}^c].
    \enalign
In the first term on the right hand side, the event $\{v \ct{A} x\}$ occurs
on $B$ but does not occur on $\gr$.  This implies that there is an occupied
path from $v$ to $x$ that exits $B$ (and does not contain a vertex in $A$),
as required.
The second term on the right hand side is zero.
To see this, we first observe that the event ${\rm NP}^c$ implies
that (on $\gr$) there is an
occupied pivotal bond $(u',v')$ for $v \conn x$
such that $v \ct{A} u'$.  The bond $(u',v')$ must also be pivotal
for the connection from $v$ to $x$ on $B$.  Moreover, since
$E'(v,x;A)$ occurs on $B$, it must be that $v$ is connected to $u'$ on $B$
and $\{v \ct{A} u'\}^c$ occurs on $B$.  This contradicts $v \ct{A} u'$
on $\gr$, and hence the second term is indeed zero.
This completes the proof of \refeq{1stdif}.

\smallskip \noindent \emph{Proof of \refeq{2nddif}.}
We begin with the identity
    \eq
    \lbeq{secterm}
    I[E^\smallsup{N}]
    -
    I[E^\smallsup{N}_\sR]
    =
    \sum_{j=1}^N \prod_{i=0}^{j-1} I[E_{i,\sR}]
    \big(I[E_{j}]-I[E_{j,\sR}]\big)\prod_{i=j+1}^N
    I[E_i],
    \en
in which the absent term with $j=0$ is equal to zero.
It suffices to show that $|I[E_{j}]-I[E_{j,\sR}]|$ is bounded above
by the indicator that $B_j$ intersects $\tilde C_{j-1} \setminus
\tilde C_{j,R}$, multiplied by either $I[E_{j}]$ or $I[E_{j,\sR}]$.
The former gives rise to the summand on the right hand side of
\refeq{2nddif}, while the latter, in combination with the products
over $i$ in \refeq{secterm}, ensures that all connections necessary
to imply $F^\smallsup{N}$ are present.  We proceed to obtain this
estimate for $|I[E_{j}]-I[E_{j,\sR}]|$.

For $A\subset A'$, we write
    \eqalign
    &I[E'(v,x;A')]-I[E'(v,x;A)]=I[E'(v,x;A)^c\cap
    E'(v,x;A')]-I[E'(v,x;A)\cap E'(v,x;A')^c],
    \lbeq{EAA'}
    \enalign
and we consider the two terms on the right hand side separately.
For the first term, we use the fact that if NP occurs for $A'$ then it
also occurs for $A$.  Therefore, the event in the first term implies that
$E'(v,x;A')\cap \{v\ct{A'}x\}\cap \{v\ct{A}x\}^c$ occurs.
In particular, there must be an occupied path from $v$ to $x$
that intersects $A'\setminus A$.
Similarly, for the second term in \refeq{EAA'},
we have $\{v \ct{A} x\} \subset \{ v \ct{A'} x\}$, so the event of the
second term implies that $E'(v,x;A)$ occurs, and that
the NP condition holds for
$A$, but not for $A'$. The latter implies that, as required,
there is an occupied
path from $v$ to $x$ containing an element in $A' \setminus A$.
This completes the proof of \refeq{2nddif}, and completes the proof of the
lemma.
\qed

\vskip 0.5cm
\noindent
{\em Proof of Lemma \ref{lem-E3}.}
We first argue that if $\vec D = \vec d$, but
$\{\vec D < \vec d\}$ on $\Vbold_{\vec d,\sM}^\smallsup{N}$,
then there must be an occupied path of length at least $10(M+1)$
in some $\Pcal_j$,
and hence,
\eq
\lbeq{lem-E3a}
    I[\{\vec D < \vec d\}\text{ on } \Vbold_{\vec d,\sM}^\smallsup{N}]
    I[\vec D = \vec d]
    \leq
    \sum_j I[\exists \text{ occupied }
    \omega \in \Pcal_j \setminus \Pcal_{j,10(M+1)}]
    I[\vec D = \vec d].
\en
Suppose, to the contrary, that
all occupied paths in each $\Pcal_j$ have length at most $10(M+1)$.
Let $R\geq r(M) =10(M+1)^2$.
Then if $\vec D = \vec d$, it must also be the case
that $\{\vec D = \vec d\}$ on
$\Vbold_{\vec d,R}^\smallsup{N}$, since the paths that determine
$\vec D$ travel at most a distance $10(M+1)$ in each of the $N+1$ expectations,
and hence travel a total distance at most $10(M+1)(N+1) \leq R$.
This proves \refeq{lem-E3a}.

We substitute \refeq{lem-E3a}
into the left hand side of \refeq{E3lem}, perform the sum
over $\vec d$ using the indicator $I[\vec D = \vec d]$, and apply
Proposition~\ref{prop-PibdsM}, to obtain the desired estimate.
\qed

\vskip 0.5cm
\noindent
{\em Proof of Lemma \ref{lem-E2}.}
We again
assume some familiarity with the methods of \cite[Section~4]{BCHSS04b}.
Fix $M$ and $N \leq M$, and fix a positive number $a<1$.

Suppose that $F^\smallsup{N}$ occurs and that $\|\vec D\| >R$.
By Proposition~\ref{prop-PibdsM}, we need only consider the case in which
all occupied paths in each $\Pcal_j$ have length at most
$10(R^a+1)$, since the complement
obeys the desired estimate with $g_R=R^a+1$.
We make this assumption throughout the proof.

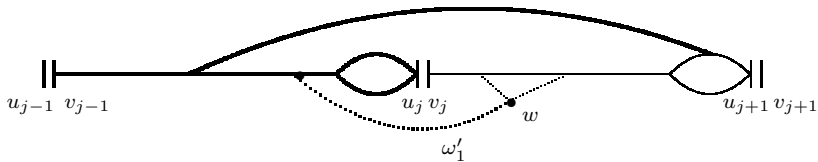
\begin{figure}[t]
\vskip2cm
\begin{center}
\setlength{\unitlength}{0.0070in}%
\begin{picture}(800,-100)
\put(580,-5){${\scriptstyle u_{j+1}}$}
\put(618,-5){${\scriptstyle v_{j+1}}$}
\put(601,10){\line(0,1){20}}
\put(609,10){\line(0,1){20}}
\Thicklines
\put(73,10){\line(0,1){20}}
\put(80,10){\line(0,1){20}}
\put(45,-5){${\scriptstyle u_{j-1}}$}
\put(88,-5){${\scriptstyle v_{j-1}}$}
\put(80,20){\line(1,0){210}}
\qbezier(290,20)(320,-10)(350,20)
\qbezier(290,20)(320,50)(350,20)
\qbezier(180,20)(375,110)(570,35)
\thicklines
\qbezier[45](260,18)(345,-50)(420,0)
\thinlines
\qbezier[15](420,0)(440,9)(460,18)
\qbezier[10](420,0)(410,9)(400,18)
\put(352,10){\line(0,1){20}}
\put(340,-5){${\scriptstyle u_{j}}$}
\put(360,-5){${\scriptstyle v_{j}}$}
\put(360,10){\line(0,1){20}}
\put(360,20){\line(1,0){180}}
\qbezier(540,20)(570,-10)(600,20)
\qbezier(540,20)(570,50)(600,20)
\put(370,-40){${\scriptstyle \omega_1'}$}
\put(260,15){${\scriptstyle \bullet}$}
\put(418,-5){${\scriptstyle \bullet}$}
\put(430,-15){${\scriptstyle w}$}
\end{picture}
\vskip1cm

\end{center}
\caption{Examples of disjoint connections required by $\omega_1'$.
The path $\omega_1'$ is a level-$j$ path, whereas the other two dotted
paths are part of the level-$(j+1)$ backbone.}
\label{fig-ii}
\end{figure}

Given a bond configuration,
we can select a sequence of occupied level-$j$ paths $\eta^j_i
\in \Pcal_j$ (a path may consist of a single vertex and the paths need
not be disjoint),
for $i = 1,2,3$ and $j=0,\ldots,N$,
which together ensure that $F^\smallsup{N}=\cap_{j=0}^N F_j$ occurs
(see \refeq{FFj}).  In Figure~\ref{fig-ii}, the paths $\eta_1^j,\eta_2^j,\eta_3^j$
are the two paths joining $v_{j-1}$ to $u_j$ and the path joining $v_{j-1}$
to $B_{j+1}$.
Denote the union of the vertices in $\eta_1^j,\eta_2^j,\eta_3^j$ by
$A_j$, and let $A=\cup_{j=0}^N A_j$.  By our assumption,
the set
$A$ explores at most $(N+1)30(R^a+1) \leq 30(M+1)(R^a+1)$ dimensions.
We consider the case where $R$ is large (depending on $M$),
 so that in particular
$R >30(M+1)(R^a+1)$.  This
implies that there must be additional occupied paths in $\cup_{j=0}^N \Pcal_j$
that explore additional dimensions.  In fact, there must be some $j$
for which the number of dimensions explored at level-$j$
exceeds $R'=(M+1)^{-1}[R-30(M+1)(R^a+1)]$,
and hence the number of these paths exceeds $R''=R'/10(R^a+1)$.
We fix such a $j$.
More precisely,
given a bond configuration for which $\|\vec D\| >R$, we can
find a $j$ and a
sequence of occupied paths $\omega_1,\ldots, \omega_{R''} \in \Pcal_j$,
such that $\omega_1$ enters a dimension not entered by $A$,
and, for $l \geq 2$,
$\omega_l$ enters a dimension not entered by $A \cup (\cup_{k<l}\omega_k)$,
where the union refers to a union of vertices.
Note that since $a<1$, $R'' \to \infty$ as $R \to \infty$ with $M$ fixed.

The paths $\eta^j_i$ ($j=0,\ldots,N$, $i=1,2,3$)
ensure that the disjoint connections required by
the event $F^\smallsup{N}$ occur.  If we were to neglect the fact
that the paths $\omega_l$ are occupied, an application of the BK
inequality would lead to the standard diagrammatic estimates for
$\hat\Pi^\smallsup{N}_p$, as described, e.g., in \cite[Section~4]{BCHSS04b}.
With this in mind,
we take the paths $\omega_l$ into account sequentially, as follows.

First, since $\omega_1$ enters a dimension not yet entered by $A_j$
(where $j$ is the special level fixed above),
there is a vertex in $\omega_1$ that is not in any of the $\eta_i^j$
($i=1,2,3$).
By following $\omega_1$ forward and backwards until it hits $A_j$
or $B_{j+1}$ for the first time, we
obtain a portion $\omega_1'$
of $\omega_1$ that begins in $A_j$ and ends in
either $A_j$ or $B_{j+1}$ and that is disjoint from the paths $\eta_i^j$.

If $\omega_1'$ both begins and ends in $A_j$, then it has the effect of
connecting two vertices on the paths $\eta_1^j,\eta_2^j,\eta_3^j$ by
an occupied path which is disjoint from these paths.  If we apply
the BK inequality in this situation, we produce Feynman diagrams
that are constructed from those bounding $\hat\Pi^\smallsup{N}_p$
by adding two vertices on diagram lines and joining them by
a new line.  In other words, we replace the product (say)
$\tau_p(y_2- y_1)\tau_p(y_4- y_3)$ by
\eq
    \sum_{a,b}
    \tau_p(y_2- a)\tau_p(a- y_1)\tau_p(y_4- b)
    \tau_p(b- y_3)\tau_p^\smallsup{1}(b-a)
\en
as in \refeq{yab}.  As in the proof of Proposition~\ref{prop-PibdsM},
we bound the factor $\sup_{a,b}\tau_p^\smallsup{1}(b-a)$ by $K\cn^{-1}$
using \refeq{tauiT},
and we are left with a diagram with two additional vertices $a,b$.
This diagram was bounded by a constant in
the proof of Proposition~\ref{prop-PibdsM}.

If $\omega_1'$ begins in $A_j$ and ends in $B_{j+1}$, then
this has the effect of augmenting $\eta_1^j,\eta_2^j,\eta_3^j$
with a disjoint path from a vertex on one of these paths to
the end of $\omega_1'$ in $B_{j+1}$.  Call this latter endpoint
$w$.   Since $w$ is in the first entry of $\omega_1'$
into $B_{j+1}$, we can augment the level-$(j+1)$
paths $\eta_1^{j+1},\eta_2^{j+1},\eta_3^{j+1}$ by a disjoint
level-$(j+1)$ path passing through $w$, as indicated in
Figure~\ref{fig-ii}.  If we apply the BK inequality to this configuration,
the result is a sum of Feynman diagrams, with extra lines due to
these additional disjoint connections.
We can begin to bound this diagram by extracting a
factor $\sup_{a,b}\tau_p^\smallsup{1}(b-a) \leq K\cn^{-1}$ from
the connection due to $\omega_1'$, and then a factor
$T_p^\smallsup{0,2} \leq K_{0,2}$ due to the level-$(j+1)$
connections that contain $w$ and are bond-disjoint from the level-$(j+1)$
connections $\eta_i^{j+1}$.
This leaves a standard $\hat\Pi^\smallsup{N}_p$ diagram
with at most three extra vertices, and this is bounded by a constant
by standard bounds, using \refeq{Tij-bd} (provided we take the
dimension sufficiently large).

The above explains the procedure if there were only one path
$\omega_1$.  However, we are interested in the situation where
there is a large number $R''$ of paths $\omega_l$.
In this case, we first find the path $\omega_1'$ as above.
Because $\omega_2$ explores a dimension not entered by
$A \cup \omega_1$, we can find a
subpath $\omega_2'$ of $\omega_2$ that starts at a vertex
in $\eta_1^j,\eta_2^j,\eta_3^j$ or $\omega_1'$ and ends either
in this set or in $B_{j+1}$, and that is disjoint from
$\eta_1^j,\eta_2^j,\eta_3^j$ and $\omega_1'$.
If $\omega_2'$ ends in $B_{j+1}$, then we can find disjoint
level-$(j+1)$ paths as explained above.
If we had only these two paths $\omega_1', \omega_2'$,
we would apply the BK inequality
as usual to produce a Feynman diagram, and then estimate this
diagram by first bounding the lines created
by $\omega_2'$, obtaining a factor $\cn^{-1}$ from the
fact that $\omega_2'$ takes at least one step, then bounding
the lines created by $\omega_1'$, obtaining a
second factor $\cn^{-1}$.  This leaves a standard $\hat\Pi^\smallsup{N}_p$
diagram with at most six
additional vertices, and this can be bounded by a constant.

The general case is handled similarly.  Each $\omega_l'$ gives
rise to a diagram line that produces
a factor $\cn^{-1}$, creating an overall factor
$\cn^{-R''}$.  In bounding diagram lines iteratively, we may
encounter lines with extra vertices (where lines already bounded
were previously attached).  However, the number of these vertices
on any one line is less than $4R''$, since each $\omega_l'$ adds
in total at most four vertices to the diagram, as in Figure~\ref{fig-ii}.
After all the lines due to the
$\omega_l'$ have been bounded using suprema, we are left with
a standard $\hat\Pi^\smallsup{N}_p$ diagram, again with
at most $4R''$ extra vertices.  This is bounded by a constant
depending on $R''$ and $M$, for $n$ sufficiently large, using
\refeq{Tij-bd}.  The number of diagrams produced depends only on $M$
and $R''$.
Since $R'' < R$, we end up with an overall bound that is a
constant multiple of $\cn^{-R''}$, where the constant
depends on $R$ and $M$.

This completes the proof.
\qed

\section*{Acknowledgements}
We thank
Christian Borgs, Jennifer Chayes and Joel Spencer for many
stimulating discussions related to this work.
The work of RvdH
was supported in part by Netherlands Organisation for
Scientific Research (NWO), and was carried out in part at
Delft University of Technology, at the University of British Columbia,
and at Microsoft Research.
The work of GS was supported in part by NSERC of Canada,
by a Senior Visiting Fellowship at the Isaac Newton
Institute funded by EPSRC Grant N09176, by EURANDOM,
and by the Thomas Stieltjes Institute.


\end{document}